	\documentclass{article}
	\usepackage{latexsym}

	\newcommand{\al}{\alpha}
	\newcommand{\del}{\delta}
	
	\newcommand{\lam}{\lambda}
	\newcommand{\vph}{\varphi}
	\newcommand{\sig}{\sigma}
	\newcommand{\tht}{\theta}
	\newcommand{\om}{\omega}
	\newcommand{\Del}{{\mathit \Delta}}
	\newcommand{\Gam}{{\mathit \Gamma}}
	\newcommand{\Lam}{{\mathit \Lambda}}
	\newcommand{\Sig}{{\mathit \Sigma}}
	\newcommand{\Om}{{\mathit \Omega}}

	\font\smbbb=msbm5 %at 11 pt
	\font\bbb=msbm7 %at 11 pt
	\font\BBB=msbm10 %at 11 pt
	\newcommand{\NN}{{\mbox{\BBB{N}}}}
	\newcommand{\ZZ}{{\mbox{\BBB{Z}}}}
	\newcommand{\QQ}{{\mbox{\BBB{Q}}}}
	
	\newcommand{\RE}{{\mbox{\BBB{R}}}}
	\newcommand{\smco}{{\mbox{\smbbb{C}}}}
	\newcommand{\co}{{\mbox{\bbb{C}}}}
	\newcommand{\CO}{{\mbox{\BBB{C}}}}
	\font\frak=eufm10 at 11 pt
	\newcommand{\g}{\mbox{\frak{g}}}
	
	\newcommand{\gt}{\mbox{\frak{t}}}
	\newcommand{\CP}{\CO P}
	\newcommand{\OO}{{\cal O}}
	\newcommand{\bfgam}[1]{\mbox{\boldmath$\mathit{\Gamma}^{{#1}}$}}

	\newcounter{sect}
	\newcommand{\sect}[1]{\vspace{4.5ex}\addtocounter{sect}{1}
		\begin{flushleft}
		{{\large\bf \arabic{sect}. {#1}}}
		\end{flushleft}
		\setcounter{thm}{0}
		\setcounter{equation}{0}
		\def\theequation{\arabic{sect}.\arabic{equation}}
		\def\thefigure{\arabic{sect}.\arabic{figure}}}
	
	\newtheorem{thm}{Theorem}[sect]
	\newtheorem{prop}[thm]{Proposition}
	\newtheorem{lemma}[thm]{Lemma}
	\newtheorem{cor}[thm]{Corollary}
	\newtheorem{defn}[thm]{Definition}
	\newtheorem{ex}[thm]{Example}
	\newtheorem{rmk}[thm]{Remark}
	\newcommand{\proof}[1]{\noindent {\em Proof.}$\quad$ {#1} $\hfill\Box$}
	\newcommand{\be}{\begin{equation}}
	\newcommand{\ee}{\end{equation}}
	\newcommand{\bea}{\begin{eqnarray}}
	\newcommand{\eea}{\end{eqnarray}}
	\newcommand{\nno}{\nonumber \\}
	\newcommand{\sep}[1]{\!\!\!\! &{#1}& \!\!\!\! }
	\newcommand{\eq}{\sep{=}}
	\newcommand{\vc}{\sep{ }}

	\newcommand{\bra}{\langle}
	\newcommand{\ket}{\rangle}
	\newcommand{\m}{\backslash}
	\newcommand{\inv}[1]{\frac{1}{#1}}
	\newcommand{\hf}{{\textstyle \inv{2}}}
	\newcommand{\e}[1]{e^{{#1}}}
        \newcommand{\ii}{\sqrt{-1}}
	\newcommand{\ch}{{\,\mathrm{char}}}
	\newcommand{\tr}{{\,\mathrm{tr}}}
	\newcommand{\ind}{{\mathrm{Ind}}}
	\newcommand{\supp}{{\,\mathrm{supp}}}
	\newcommand{\im}{{\mathrm{Im}\,}}
	\newcommand{\set}[2]{\{{#1}\,|\,{#2}\}}
	\newcommand{\dr}{d}
	\newcommand{\pdb}{\bar{\partial}}

	\newcommand{\wb}{\bar{w}}
	\newcommand{\pC}{^{p,C}}
	\newcommand{\sC}{^{\sig,C}}
	\newcommand{\pk}{^p_k}
        \newcommand{\lpk}{{\lam^p_k}}
        \newcommand{\lSk}{{\lam^S_k}}
	
	\newcommand{\SL}{^S_\Lam}
	\newcommand{\kL}{^k_\Lam}
	\newcommand{\wl}{w\lam^S_k}
	\newcommand{\wL}{w(\Lam+\rho)-\rho}
	\newcommand{\nwl}{{n^C(\{\wl\})}}
	\newcommand{\prwl}{\prod_{\wl\in C^*}\inv{1-\e{-\wl}}
			\prod_{\wl\in-C^*}\frac{\e{\wl}}{1-\e{\wl}}}
	\newcommand{\TC}{T^\co}
	
	\newcommand{\HC}{H^\co}
	\newcommand{\Gc}{G^\smco}
	\newcommand{\GC}{G^\co}
	\newcommand{\Dp}{\Del^+\m\Del^+}

	\newcommand{\dom}{{\cal D}}
	\newcommand{\lat}{{\cal L}}
	\newcommand{\ring}{\ZZ[\lat^*]}
	\newcommand{\sumk}{\sum^n_{k=0}}
        \newcommand{\prodk}{\prod^n_{k=1}}
	\newcommand{\sump}{\sum_{p\in F}}
	\newcommand{\sumw}{\sum_{w\in W}}
        \newcommand{\hme}{{H^k(M,\OO(E))}} 
	\newcommand{\ka}{K\"ahler }
	\newcommand{\pra}[1]{\prod_{\al\in\Del^+_{{#1}}}(1-\e{-\al})}

\begin{document}
$\!\,{}$

        \vspace{-5ex}

        \begin{flushright}
{\tt dg-ga/9602008} (February, 1996)\\
MSRI preprint No.\ 1996-013\\
Revised (December, 1997)
        \end{flushright}

	\begin{center}
	{\LARGE\bf Equivariant Holomorphic Morse Inequalities II:\\
	\vspace{1ex}
		Torus and Non-Abelian Group Actions}\\

	\vspace{4ex}
	{\large\rm Siye Wu}\footnote{Current address: Department of 
	Pure Mathematics, University of Adelaide, Adelaide 5005, Australia.
	E-mail address: {\tt swu@maths.adelaide.edu.au}}

	{\em Mathematics Section, International Centre for Theoretical Physics,
	Trieste I-34100, Italy\\
	and\\
	Mathematical Sciences 	Research Institute, 1000 Centennial Drive,
	Berkeley, CA 94720, USA}
	\end{center}

	\vspace{3ex}

        \begin{quote}
{\small {\bf Abstract.}
We extend the equivariant holomorphic Morse inequalities of circle actions
to cases with torus and non-Abelian group actions on holomorphic vector
bundles over \ka manifolds and show the necessity of the \ka condition.
For torus actions, there is a set of inequalities for each choice of action
chambers specifying directions in the Lie algebra of the torus.
We apply the results to invariant line bundles over toric manifolds.
If the group is non-Abelian, there is in addition an action of the Weyl group
on the fixed-point set of its maximal torus.
The sum over the fixed points can be rearranged into sums over the
Weyl group (having incorporated the character of the isotropy representation
on the fiber) and over its orbits.}
        \end{quote}

        \vspace{3ex}

\sect{Introduction}

Index theorems express analytical indices of elliptic complexes in terms of
topological invariants;
information on the individual cohomology groups are usually obtained
with the aid of vanishing theorems.
Taking the de Rham complex for example, the Euler number is not enough to
determine the Betti numbers.
However, if we consider a Morse function, then the Morse inequalities bound
each Betti number by the data of the critical points.
In this paper, we consider a holomorphic setting in which a compact group
acts holomorphically on a holomorphic vector bundle over a \ka manifold.
The index theorem is the Atiyah-Bott fixed-point formula [\ref{AB}]
(when the fixed points are isolated),
which expresses the equivariant index, the alternating sum of the characters
of the Dolbeault cohomology groups, in terms of the fixed-point data.
The corresponding equivariant holomorphic Morse inequalities when the group
is the circle group was obtained by Witten [\ref{W}] and was first proved 
analytically using the heat kernel method by Mathai and the present author
[\ref{MW}] when the fixed points are isolated.
In this paper, we extend the result to cases with torus and non-Abelian
group actions.
We also show that the \ka condition is essential for such Morse-type
inequalities although, in contrast, not necessary for the equivariant
index formula of Atiyah and Bott.

Let $M$ be a compact \ka manifold of complex dimension $n$ and $E$,
a holomorphic vector bundle over $M$.
The Dolbeault cohomology groups $H^*(M,\OO(E))$ with coefficients in $E$ are
cohomologies of the twisted Dolbeault complex $(\Om^{0,*}(M,E),\pdb_E)$
and are independent of the choice of holomorphic connections.
Let $G$ be a compact, connected Lie group whose Lie algebra is
denoted by ${\mathrm{Lie}}(G)$. 
Let $\g=\ii\,{\mathrm{Lie}}(G)$.
We assume that $G$ acts holomorphically and effectively on $M$
preserving the \ka form $\om$.
If the fixed-point set of a maximal torus in $G$ is non-empty, then
the $G$-action is Hamiltonian [\ref{Fr}], i.e.,
there is a moment map $\mu\colon M\to\g^*$ such that for any $x\in\g$,
the corresponding vector field $V_x$ on $M$ satisfies
$i_{V_x}\om=\dr\bra\mu,x\ket$.
If the action of $g\in G$ on $M$ (still denoted by $g$) can be lifted
holomorphically to $\tilde{g}$ on $E$, then $G$ acts on the space of 
sections $\Gam(M,E)$ by $g\colon s\mapsto\tilde{g}\circ s\circ g^{-1}$
($g\in G$) and similarly on $\Om^*(M,E)$.
In this case, $G$ commutes with the twisted Dolbeault operator $\pdb_E$
and hence acts on the cohomology groups $H^k=H^k(M,\OO(E))$ ($0\le k\le n$).
The purpose of this paper is to study the decomposition of
each $H^k$ in terms of the irreducible representations of $G$.

In sections 2 and 3, the group acting on $M$ is a torus $T$.
Section 2 shows that there is a set of Morse-type inequalities
for each choice of action chambers specifying the directions
in the Lie algebra $\gt$.
This is obtained by applying the result of [\ref{W}, \ref{MW}] to various
circle subgroups of $T$.
Section 3 applies the result of the previous section to $T$-invariant line
bundles over toric manifolds (including projective spaces and Hirzebruch
surfaces).
In section 4, we demonstrate that for non-\ka manifolds,
the strong equivariant holomorphic Morse inequalities need not hold.
Violation of strong inequalities is also used to show the non-existence
of invariant \ka structures.
In section 5, the group $G$ is a general compact non-Abelian group.
The main result is obtained by applying that of section 2 to
a maximal torus $T$ of $G$, assuming that the $T$-fixed-point set $F$ is 
discrete.
The novelties are that the cohomology groups $H^k$ ($0\le k\le n$) are
representations of the non-Abelian group $G$ (hence the structure of $H^k$
is more rigid) and that there exists an action of the Weyl group $W$ on $F$
(hence the set $F$ can be organized into $W$-orbits).

Throughout this paper, $\CO$, $\RE$, $\QQ$, $\ZZ$, $\NN$ will denote the
sets of complex numbers, real numbers, rational numbers, integers,
non-negative integers, respectively.

\sect{Equivariant holomorphic Morse inequalities with torus actions}

In this section, we assume that the compact Lie group that acts on the
\ka manifold is a torus group $T$.
Let $\gt=\ii\,{\mathrm{Lie}}(T)$ be the Lie algebra of $T$ and $\lat$,
the integral lattice in $\gt$;
the dual lattice $\lat^*$ in $\gt^*$ is the weight lattice.

\begin{defn}
Let $\ring$ be the formal character ring of $T$ consisting of elements 
$q=\sum_{\xi\in\lat^*}q_\xi\e{\xi}$ ($q_\xi\in\ZZ$).
We say $q\ge0$ if $q_\xi\ge0$ for all $\xi\in\lat^*$.
The support of $q$ is the set $\supp(q)=\set{\xi\in\lat^*}{q_\xi\ne0}$.
Let $Q(t)=\sumk q_kt^k\in\ring[t]$ be a polynomial of degree $n$ with
coefficients in $\ring$.
We say $Q(t)\ge0$ if $q_k\ge0$ in $\ring$ for all $k$.
For two such polynomials $P(t)$ and $Q(t)$, we say $P(t)\le Q(t)$ if
$Q(t)-P(t)\ge0$.
\end{defn}

For example, if $V$ is a finite dimensional representation of $T$, then
its character $\ch(V)\ge0$ in $\ring$.
Let the {\em support} of $V$, denoted by $\supp(V)$, be the set
$\supp(\!\ch(V))$ of weights whose multiplicity in $V$ is non-zero.
For any $\tht\in\gt$, there is a homomorphism $\ring\to\CO$
given by $\e{\xi}\mapsto\e{\ii\bra\xi,\tht\ket}$.
For instance, $\ch(V)\mapsto\tr_V\e{\ii\tht}$ under this homomorphism.
Another important type of elements in $\ring$ is given by the series
	\be\label{series}
\frac{\e{\eta}}{1-\e{\xi}}\stackrel{{\rm def.}}{=}
\sum_{k=0}^\infty\e{k\xi+\eta},
\quad\xi,\eta\in\lat^*.
	\ee
We emphasize here that in (\ref{series}) the left-hand side is a notation
for the formal series on the right-hand side.

Recall we assumed that $T$ acts holomorphically and effectively on a 
compact \ka manifold $M$ with a non-empty and discrete fixed-point set $F$.
We also assume that the $T$-action preserves the \ka form $\om$ and hence is
Hamiltonian [\ref{Fr}]; let $\mu\colon M\to\gt^*$ be the moment map.
For any $p\in F$, let $\lam^p_1,\cdots,\lam^p_n\in\lat^*\backslash\{0\}$
be the weights of the isotropy representation of $T$ on $T_pM$
(Figure \ref{definition}(a)).
The hyperplanes $(\lam^p_k)^\perp\subset\gt$ cut $\gt$ into open polyhedral
cones called {\em action chambers} [\ref{GLS}, \ref{GP}, \ref{PW}].
(When $M$ is a coadjoint orbit, the action chambers are precisely the
Weyl chambers.)
We fix a {\em positive} action chamber $C$ (Figure \ref{definition}(b)).
Let $\lam\pC_k=\pm\lam^p_k$ be the {\em polarized weights}, with the sign 
chosen so that $\lam\pC\in C^*$.
(Here $C^*$ is the dual cone in $\gt^*$ defined by 
$C^*=\set{\xi\in\gt^*}{\bra\xi,C\ket>0}$.)
We define the {\em polarizing index} $n\pC$ of $p$ with respect to $C$
as the number of weights $\lam^p_k\in-C^*$.
We also assume that there is a holomorphic vector bundle $E$ over $M$
on which the $T$-action lifts holomorphically.
The torus $T$ acts on the fiber $E_p$ over $p\in F$ and
on the cohomology groups $H^k=\hme$ ($0\le k\le n$).
Following [\ref{W}, \ref{MW}], we denote $E_p(\tht)=\tr_{E_p}\e{\ii\tht}$
($p\in F$) and $H^k(\tht)=\tr_{H^k}\e{\ii\tht}$.

\begin{thm}\label{TORUS}                                                
For each choice of the positive action chamber $C$,
we have the strong equivariant holomorphic Morse inequalities
        \be\label{strong}
\sump t^{n\pC}\ch(E_p)\prod_{\lam\pk\in C^*}\inv{1-\e{-\lam\pk}} 
\prod_{\lam\pk\in-C^*}\frac{\e{-\lam\pC_k}}{1-\e{-\lam\pC_k}}
=\sum_k\,t^k\ch(H^k)+(1+t)Q^C(t)
        \ee                                                       
for some $Q^C(t)\ge0$ in $\ring$.
\end{thm}

\begin{rmk}
{\em The strong inequalities (\ref{strong}) imply the following
weak equivariant holomorphic Morse inequalities: 
	\be\label{weak}
\ch(H^k)\le\sum_{n\pC=k}\ch(E_p)\prod_{\lam^p_j\in C^*}\inv{1-\e{-\lam^p_j}} 
\prod_{\lam^p_j\in-C^*}\frac{\e{-\lam\pC_j}}{1-\e{-\lam\pC_j}}
\quad(0\le k\le n).
	\ee
Setting $t=-1$ in (\ref{strong}), we recover the Atiyah-Bott fixed-point 
theorem
        \be\label{ab}                                              
\sump\frac{\ch(E_p)}{\prodk(1-\e{-\lpk})}=\sumk(-1)^k\ch(H^k)=\ind.
        \ee                                                                    
}
\end{rmk}

\begin{rmk}
{\em If $T$ is the circle group $S^1$, then 
$\gt=\ii\,{\mathrm{Lie}}(S^1)\cong\RE$.
There are only two action chambers $C^\pm=\set{\tht\in\RE}{\pm\tht>0}$.
The polarizing index $n^{p,C^+}$ of $p$ with respect to $C^+$ is the number
of weights $\lam\pk<0$, denoted by $n_p$.
Theorem~\ref{TORUS} reduces to the results in [\ref{W}, \ref{MW}].
In particular, for $C=C^+$, (\ref{strong}) becomes
        \be\label{strong+}
\sump t^{n_p}E_p(\tht)\prod_{\lpk>0}\inv{1-\e{-\ii\lpk\tht}}\prod_{\lpk<0}
\frac{\e{-\ii|\lpk|\tht}}{1-\e{-\ii|\lpk|\tht}}
=\sumk t^kH^k(\tht)+(1+t)Q^+(\tht,t), 
        \ee
where $Q^+(\tht,t)\ge0$ in $\RE((\e{\ii\tht}))$, the ring of formal characters
of $S^1$ [\ref{MW}].}
\end{rmk}

\noindent {\em Proof of Theorem~\ref{TORUS}.}$\quad$
We notice that (\ref{strong}) is equivalent to
	\be\label{anal}
\sump t^{n\pC}E_p(\tht)
\prod_{\lam\pk\in C^*}\inv{1-\e{-\ii\bra\lam\pk,\tht\ket}} 
\prod_{\lam\pk\in-C^*}\frac{\e{-\ii\bra\lam\pC_k,\tht\ket}}                   
{1-\e{-\ii\bra\lam\pC_k,\tht\ket}}
=\sum_k\,t^kH^k(\tht)+(1+t)Q^C(t)(\ii\tht)    
        \ee
regarded as an equality of analytic functions in $\tht\in\gt^\co$ for
$\im\tht\in-C$.
(\ref{strong}) implies (\ref{anal}) because for any $\del\in C$ 
the Taylor expansions on the left hand side of (\ref{anal}) are uniformly
convergent in $\tht\in\gt^\co$ if $\im\tht\in-(\del+\overline{C})$.
Conversely, if (\ref{anal}) is true for complex $\tht\in\gt^\co$ with  
$\im\tht\in-C$, we take $\im\tht\to 0$ within the cone $-C$.
The equality (\ref{strong}) of formal series follows from the uniqueness
of the Fourier expansion of tempered distributions [\ref{Ho}].
For $S^1$-actions [\ref{W}, \ref{MW}], (\ref{strong+}) is a true equality
both in $\RE((\e{\ii\tht}))$ and as analytic functions for $\im\tht<0$.
The rest of the proof, which is similar to that of [\ref{PW}, Theorem~2.2],
shows (\ref{anal}) using (\ref{strong+}).
By analyticity, it suffices to prove (\ref{anal}) for purely imaginary $\tht$. 
Pick an integral point $\tht_1\in C\cap\lat$. 
Since $\bra\lam^p_k,\tht_1\ket\ne0$ for all $p\in F$ and $1\le k\le n$,
$h=\bra\mu,\tht_1\ket$ generates a Hamiltonian $S^1$-action on $M$ 
with the same fixed-point set $F$.
Moreover the weights of $S^1$ at $p\in F$ are $\bra\lam^p_k,\tht_1\ket$
($1\le k\le n$) and $|\bra\lam^p_k,\tht_1\ket|=\bra\lam\pC_k,\tht_1\ket>0$
for all $p\in F$ and $1\le k\le n$.
(\ref{strong+}) implies (\ref{anal}) for $\tht=-\ii s\tht_1$ ($s>0$).
The result follows from continuity since such $\tht$'s form a dense subset 
of $-\ii C$.							$\hfill\Box$

The inequalities (\ref{weak}) or (\ref{strong}) show that the multiplicities
of weights in the cohomology groups $H^k$ ($0\le k\le n$) are constrained
by the fixed-point data. 
Given an action chamber $C$, the support of $H^k$ is contained in 
a suitably shifted cone $-C^*$ in $\lat^*$.
By choosing different chambers, it is possible to bound $\supp H^*$
($0\le k\le n$) in various directions in $\gt^*$.
We need the following definition. (See Figure \ref{definition}(c),(d).)

\begin{defn}\label{region}
For a given choice of positive action chamber $C$ and $p\in F$, let
	\be
\Gam\pC=\set{\xi-\textstyle{\sum_{k=1}^n}r_k\lam\pC_k}{\xi\in\supp(E_p),
r_k\ge0\mbox{ {\em and} }r_k>0\mbox{ {\em if} }\lam\pk\in-C^*}
	\ee
and
	\be
\Gam^{k,C}=\bigcup_{p\in F,n\pC=k}\Gam\pC.
	\ee
We set $\Gam^k=\bigcap_C\Gam^{k,C}$.
\end{defn}

\begin{prop}\label{SUPP}
For any $0\le k\le n$,\\
1.	\be
\supp(H^k)\subset\lat^*\cap\bigcap_C\Gam^{k,C}=\lat^*\cap\Gam^k;
	\ee
2.	\be
\supp(H^k)\supset\lat^*\cap\bigcup_C\Gam^{k,C}\m(\Gam^{k-1,C}\cup\Gam^{k+1,C}).
	\ee
\end{prop}

\proof{Part 1 follows from the weak inequalities with all choices of $C$.
Part 2 follows from the strong inequalities:
If $\xi\in\lat^*\cap\Gam^{k,C}$ but $\xi\not\in\Gam^{k\pm1,C}$,
then the polynomial $Q^C(t)=\sumk\sum_{\xi\in\lat^*}q^k_\xi\e{\xi}t^k$
in (\ref{strong}) has coefficients $q^{k-1}_\xi=q^k_\xi=0$.
This means that $(1+t)Q^C(t)$ does not contain the term $\e{\xi}t^k$.
Hence $\xi\in\supp(H^k)$.}

Recall that a \ka manifold $(M,\om)$ is {\em quantizable} 
if $\frac{\om}{2\pi}$ represents an integral de Rham class.
In this case, a {\em pre-quantum line bundle} over $M$ is 
a line bundle whose curvature is $\frac{\om}{\ii}$.

\begin{cor}\label{CONVEX}
If $L$ is a $T$-invariant pre-quantum line bundle of a quantizable
\ka manifold $(M,\om)$, then $\supp H^k(M,\OO(L))$ is contained
in the moment polytope $\Del$ for any $k=0,\cdots,n$.
\end{cor}

\proof{Let $\mu\colon M\to\gt^*$ be the moment map.
The image $\Del=\mu(M)$ is a convex polytope [\ref{A}, \ref{GS}].
For any facet of $\Del$, there is $x\in\gt$ such that the hyperplane
$x^\perp$ contains the facet and $\bra\Del,x\ket\le0$.
Choose an action chamber $C$ with an edge containing $x$.
For any fixed point $p\in F$, the weight of the torus action on the fiber 
of the pre-quantum line bundle over $p$ is $\mu(p)\in\Del$, hence
$\bra\mu(p),x\ket\le0$.
Since all $\bra\lam\pC_k,x\ket\ge0$ we conclude that $\bra\Gam\pC,x\ket\le0$
for all $p\in F$.
By part 1 of Proposition~\ref{SUPP}, $\bra\supp(H^k),x\ket\le0$ for
any facet of $\Del$.
Since $\Del$ is convex, we get $\supp(H^k)\subset\Del$ for any $k$.}
\vspace{2ex}

In sections 3 and 5, we will apply Theorem~\ref{TORUS} and 
Proposition~\ref{SUPP} to toric manifolds and to cases with 
general non-Abelian group actions.

\sect{Applications to toric manifolds}

The toric manifolds we consider here are smooth complex manifolds $M$,
each equipped with an effective action of $T=T^n$ ($n=\dim_\co M$).
Such a manifold can be characterized combinatorially by a fan $\Sig$,
a collection of cones in $\gt=\ii\,{\mathrm{Lie}}(T)$ that satisfy
certain compatibility and integrality conditions.
(For reviews, see for example [\ref{Da}, \ref{O}, \ref{Au}, \ref{F}].)
$M$ is compact if and only if the corresponding fan $\Sig$ is complete,
i.e., the union of all the cones $|\Sig|=\gt$.
Top dimensional cones $\sig$ in the $n$-skeleton $\Sig^{(n)}$
are in one-to-one correspondent with the $T$-fixed points $p_\sig\in F$.
The weights of isotropic representations on $T_{p_\sig}M$
span the cone $-\sig^*\subset\gt^*$.
The hyperplanes containing $(n-1)$-dimensional cones cut $\gt$ into action 
chambers.
The $T$-equivariant holomorphic line bundles over $M$ are characterized by
(continuous) piecewise linear functions $\vph$ on $\Sig$ modulo $\gt^*$
(the space of globally linear functionals).
For each $\sig\in\Sig^{(n)}$, $\vph|_\sig\in\gt^*$ is the weight of 
the $T$-action on the fiber over $p_\sig$.
We denote this line bundle by $L_\vph$.
The cohomology groups $H^k=H^k(M,\OO(L_\vph))$ ($k=0,\cdots,n$) are
representations of $T$.
For each weight $\xi\in\lat^*$, the multiplicity of $\xi$ in $H^k$ can be
computed using $\xi$, $\vph$ and $\Sig$ (see Remark~\ref{dema} below).
The purpose of this section is to find information on such multiplicities
using the results of the previous section.
This method provides much geometric insight, especially on how
the multiplicity varies according to the position of the weight relative
to the image of the (generalized) moment map.
The multiplicity in the equivariant index was studied in [\ref{KT}, \ref{GK}]
with similar considerations by using the Atiyah-Bott fixed-point formula
(\ref{ab}).

We consider a toric manifold $M$ with a $T$-invariant \ka structure.
This is equivalent to the existence of a $T$-invariant ample line bundle
over $M$; the \ka structure is then induced by the projective embedding.
Such a line bundle corresponds to a strictly convex piecewise linear function.
{}From the symplectic point of view, if $M$ admits a (non-degenerate)
symplectic form, then the image of the moment map is a convex polytope
which does define a strictly convex piecewise linear function on its fan
(with a possible perturbation of the \ka form).
Therefore a symplectic toric manifold is K\"ahler.
An explicit construction of a \ka structure
using the Delzant construction [\ref{Del}] was given by [\ref{G}].

\begin{prop}\label{TORIC} Let $M$ be a compact smooth toric manifold of
complex dimension $n$ with a $T$-invariant \ka structure
and let $H^k=H^k(M,\OO(L_\vph))$ ($k=0,\cdots,n$), where the line bundle
is determined by a piecewise linear function $\vph$ on the fan. Then\\
1. If $0<k<n$, the sets $\supp(H^0)$, $\supp(H^k)$ and $\supp(H^n)$
are mutually disjoint;\\
2. $\supp(H^0)=\lat^*\cap\Gam^0$, $\supp(H^n)=\lat^*\cap\Gam^n$,
where each weight is of multiplicity $1$.
\end{prop}

\proof{Choose any action chamber $C$. 
Let $\sig_C\in\Sig^{(n)}$ be the unique cone containing $-C$.
It is clear that the polarizing indices $n\sC$ of $p_\sig\in F$ are 
given by $n^{\sig_C,C}=0$ and $n\sC>0$ for all $\sig\ne\sig_C$.
Therefore in the weak inequality (\ref{weak}) for $k=0$,
there is only one term that bounds $\ch(H^0)$.
Since $\dim_\co M=\dim\gt$, the number of isotropy weights at $p_\sig$
is equal to $n$, hence
the multiplicity of any weight in $H^0$ is not greater than $1$.
Let $C$ run through all action chambers.
Since $\Gam^{\sig_C,C}=\vph|_{\sig_C}+\sig^*_C$, we obtain
$\Gam^0=\bigcap_{\sig\in\Sig^{(n)}}(\vph|_\sig+\sig^*)$.
Using again $\dim_\co M=\dim\gt=n$, for any $\sig\in\Sig^{(n)}$
and any action chamber $C$ such that $\sig_C\ne\sig$, we have 
$(\vph|_\sig+\sig^*)\cap\Gam^{\sig,C}=\emptyset$.
Therefore $\Gam^0\cap\Gam^{\sig,C}=\emptyset$.
Taking the union of $\sig\in\Sig^{(n)}$ with $n^{\sig,C}=k>0$, we get
$\Gam^0\cap\Gam^{k,C}=\emptyset$, hence $\Gam^0\cap\Gam^k=\emptyset$ for $k>0$.
By part 1 of Proposition~\ref{SUPP}, $\supp(H^0)\cap\supp(H^k)=\emptyset$
for $k>0$.
Also, since $\Gam^0\cap\Gam^{k,C}=\emptyset$ for any $k>0$ and any $C$,
by the Atiyah-Bott fixed-point theorem (\ref{ab}) or by part 2 of
Proposition~\ref{SUPP}, we have $\supp(H^0)\supset\lat^*\cap\Gam^0$.
Hence $\supp(H^0)=\lat^*\cap\Gam^0$, each weight with multiplicity $1$.
The results on $H^n$ can be proved similarly.}

\begin{cor}\label{AMPLE}
If $\vph$ is a strictly convex piecewise linear function, then
$H^k=0$ for $k>0$ and $\ind=\ch(H^0)$;
If $-\vph$ is convex, then $H^k=0$ for $k<n$ and $\ind=(-1)^n\ch(H^n)$.
\end{cor}

\proof{In the first case, there is a $T$-invariant \ka form
such that $L_\vph$ is the pre-quantum line bundle and 
$\mu(p)=\vph|_{p_\sig}$ for all $\sig\in\Sig^{(n)}$.
The moment polytope
$\Del=\bigcap_{\sig\in\Sig^{(n)}}(\mu(p_\sig)+\sig^*)=\Gam^0$.
{}From Corollary \ref{CONVEX} and Proposition \ref{TORIC}, we obtain
$\Gam^k\subset\Del$ for all $k$, and $\Gam^0\cap\Gam^k=\emptyset$ for $k>0$.
Therefore $H^k=0$ for $k>0$, and $\ind=\ch(H^0)$, 
with $\supp(H^0)=\lat^*\cap\Del$.
The second case can be proved similarly.}

\begin{ex}
{\em Consider the projective space $\CP^n=(\CO^{n+1}\m\{0\})/\CO^*$ with 
the standard $T^n$-action given by $(u_1,\cdots,u_n)\colon
[z_1,\cdots,z_n.z_{n+1}]\mapsto[u_1z_1,\cdots.u_nz_n,z_{n+1}]$
in homogeneous coordinates.
$\gt=\ii\,{\mathrm{Lie}}(T^n)\cong\RE^n$ is spanned by $e_1,\cdots,e_n$.\\
1. The case $n=1$ has been considered in [\ref{W}] using the inequalities in
(\ref{strong+}) and their counterparts for $C=C^-$.
If $L$ is a line bundle over $\CO P^1$ with $c_1(L)=r$,
then $\dim H^0(\CO P^1,\OO(L))=r+1$, $H^1=0$ if $r\ge0$,
and $H^0=0$, $\dim H^1=|r|-1$ if $r<0$. \\
2. It is instructive to work out the case $n=2$ directly using the results
of the previous section.
There are three fixed points $p_1=[1,0,0]$, $p_2=[0,1,0]$, $p_3=[0,0,1]$
with isotropy weights $\{e^*_1,e^*_1-e^*_2\}$, $\{e^*_2,e^*_2-e^*_1\}$,
$\{-e^*_1,-e^*_2\}$, respectively.
There are six action chambers. 
Let $C$ be a chamber spanned by  $\{e_1,e_1+e_2\}$ 
(Figure~\ref{projective}(a)).
The line bundle $L=(\CO^3\m\{0\}\times\CO)/\CO^*$, where the quotient is
$(z_1,z_2,z_3,w)\sim(uz_1,uz_2,uz_3,u^rw)$, $u\in\CO^*$, has $c_1(L)=r\in\ZZ$.
The $T^2$-action lifts to $L$ by 
$(u_1,u_2)\colon[z_1,z_2,z_3,w]\mapsto[u_1z_1,u_2z_2,z_3,w]$.
The weights on $L_{p_1}$, $L_{p_2}$, $L_{p_3}$ are $\xi_1=re^*_1$,
$\xi_2=re^*_2$, $\xi_3=0$, respectively (Figure~\ref{projective}(b)).
It is easy to see that $n^{p,C}=0,1,2$ for $p=p_1,p_2,p_3$, respectively
and that $\Gam^{p_1,C}=\set{x_1e^*_1+x_2e^*_2}{x_2\ge0,x_1+x_2\le r}$,
$\Gam^{p_2,C}=\{x_1<0,x_1+x_2\le r\}$, $\Gam^{p_3,C}=\{x_1,x_2<0\}$.
On the other hand, $n^{p,-C}=2,1,0$ for $p=p_1,p_2,p_3$, respectively
and $\Gam^{p_1,-C}=\{x_1+x_2>r,x_2<0\}$,
$\Gam^{p_2,-C}=\{x_1\ge0,x_1+x_2>r\}$, $\Gam^{p_3,-C}=\{x_1,x_2\ge0\}$.
If $r\ge0$, by part 1 of Proposition~\ref{SUPP},
$\Gam^0\subset\Gam^{p_1,C}\cap\Gam^{p_3,-C}=\{x_1,x_2\ge0,x_1+x_2\le r\}$
(Figure~\ref{projective}(c)),
$\Gam^1\subset\Gam^{p_2,C}\cap\Gam^{p_2,-C}=\emptyset$
(Figure~\ref{projective}(d)),
and $\Gam^2\subset\Gam^{p_1,C}\cap\Gam^{p_3,-C}=\emptyset$ 
(Figure~\ref{projective}(e)).
By part 2 of Proposition~\ref{SUPP},
$\supp(H^0)\supset\ZZ^2\cap(\Gam^{p_1,C}\m\Gam^{p_2,C})$.
Hence $\supp(H^0)=\ZZ^2\cap\Gam^0$, where 
$\Gam^0=\{x_1,x_2\ge0,x_1+x_2\le r\}$, and $H^1=H^2=0$.
Similarly, if $r<0$, then $H^0=H^1=0$, and $H^2=\ZZ^2\cap\Gam^2$,
where $\Gam^2=\{x_1,x_2<0,x_1+x_2>r\}$.\\
3. For general $n$, $\CO P^n$ as a toric manifold can be described by a fan
$\Sig$: any $k$ vectors in $\{v_0=-\sum_{i=1}^ne_i,v_i=e_i\,(1\le i\le n)\}$
span a $k$-dimensional cone in the fan $\Sig$.
For any piecewise linear function $\vph$ on $\Sig$,
either $\vph$ or $-\vph$ is convex.
In fact by shifting a linear functional, $\vph$ can be brought to
the standard form $\vph(v_0)=-r$, $\vph(v_i)=0$ ($1\le i\le n$);
$r$ is the first Chern number of the line bundle.
$\vph$ is convex if and only if $r\ge 0$.
In this case, Corollary~\ref{AMPLE} implies that $H^k=0$ for $k>0$
and that $\supp(H^0)$ contains the integer points in the simplex
$\set{\sum_{i=1}^nx_ie_i^*}{x_i\ge0\,(1\le i\le n),\sum_{i=1}^nx_i\le r}$
in $\gt^*\cong\RE^n$, each with multiplicity $1$.
If $r<0$, then $H^k=0$ for $k<n$, and $\supp(H^n)=\ZZ^2\cap
\set{\sum_{i=1}^nx_ie_i^*}{x_i<0\,(1\le i\le n),\sum_{i=1}^nx_i>r}$,
with multiplicity $1$.}
\end{ex}

\begin{cor}\label{2D}
If $\dim_\co M=2$, then
	\be
\supp(\ind)=\supp(H^0)\cup\supp(H^1)\cup\supp(H^2)
	\ee
is a disjoint union.
\end{cor}

\proof{This follows from part 1 of Proposition~\ref{TORIC}
since the only choice of $k$ is $k=1$.}

\begin{ex}
{\em We consider the Hirzebruch surface, given by a fan $\Sig$ in $\RE^2$
whose $1$-skeleton $\Sig^{(1)}$ is spanned by vectors
$v_1=e_1$, $v_2=e_2$, $v_3=-e_1$, $v_4=-ae_1-e_2$ ($a\in\NN$).
We choose an action chamber $C$ spanned by $\{-e_1,-ae_1-e_2\}$
(Figure~\ref{hirzebruch}(a)).
Consider a piecewise linear function given by 
$\vph(v_1)=\vph(v_2)=0$, $\vph(v_3)=-r$ and $\vph(v_4)=-s$.
We assume that $r,s>0$.
It is easy to see that if $s\ge ar$, then $\vph$ is convex, hence
$H^1=H^2=0$, $\supp(H^0)=\ZZ^2\cap\Gam^0$, with multiplicity $1$,
where $\Gam^0=\{0\le x_1\le s, 0\le x_2\le r-ax_1\}$ is a $4$-gon.
If $s<ar$, since isotropy weights of the four fixed points $p_i$ 
($1\le i\le 4$) are the negatives of $\{e^*_1,e^*_2\}$, $\{-e^*_1,e^*_2\}$,
$\{-e^*_2,-e^*_1+ae^*_2\}$, $\{-e^*_2,e^*_1-ae^*_2\}$, respectively
(Figure~\ref{hirzebruch}(b)), we have $n^{p,C}=0,1,2,1$ for 
$p=p_1,p_2,p_3,p_4$, and
$\Gam^{p_1,C}=\{x_1,x_2\ge0\}$, $\Gam^{p_2,C}=\{x_1>r,x_2\ge0\}$,
$\Gam^{p_3,C}=\{x_1>r,ax_1+x_2>s\}$, $\Gam^{p_4,C}=\{x_1\ge0,ax_1+x_2>s\}$.
On the other hand, $n^{p,-C}=2,1,0,1$ for $p=p_1,p_2,p_3,p_4$, respectively,
and $\Gam^{p_1,-C}=\{x_1,x_2<0\}$, $\Gam^{p_2,-C}=\{x_1\le r,x_2<0\}$,
$\Gam^{p_3,-C}=\{x_1\le r,ax_1+x_2\le s\}$,
$\Gam^{p_4,-C}=\{x_1<0,ax_1+x_2\le s\}$.
Using Proposition~\ref{SUPP}, we obtain $\supp(H^0)=\ZZ^2\cap\Gam^0$,
$\supp(H^1)=\ZZ^2\cap\Gam^1$ (the equalities follow from part 2), where 
$\Gam^0=\Gam^{p_1,C}\cap\Gam^{p_3,-C}=\{x_1,x_2\ge0, ax_1+x_2\le r\}$
(Figure~\ref{hirzebruch}(c)),
$\Gam^1=(\Gam^{p_2,C}\cup\Gam^{p_4,C})\cap(\Gam^{p_2,-C}\cup\Gam^{p_4,-C})
=\{x_1\le r,x_2<0, ax_1+x_2>s\}$ (Figure~\ref{hirzebruch}(d)),
and $H^2=0$ since $\Gam^2\subset\Gam^{p_3,C}\cap\Gam^{p_1,-C}=\emptyset$
(Figure~\ref{hirzebruch}(e)).
The results are in accord with Corollary~\ref{2D}.
An alternative way of obtaining the same results is to calculate the 
equivariant index (\ref{ab}) and use Corollary~\ref{2D}.}
\end{ex}

\begin{rmk}
{\em We make several observations on the multiplicity of weights in $H^k$
when $0<k<n$.\\
1. First, the multiplicity of a weight in $H^k$ ($0<k<n$) is 
not necessarily $1$, even when $\dim_\co M=2$.
For example, it was shown in [\ref{KT}, Example~5.6] that for a certain line
bundle over the blow-up of $\CP^2$ at three points, a weight can appear
in the index with multiplicity $-2$.
Using Corollary~\ref{2D}, one concludes that the multiplicity of that weight
in $H^1$ is $2$.\\
\noindent
2. When $\dim_\co M\ge3$, a weight can appear simultaneously
in cohomology groups of different degrees.
To see this, we construct a fan $\Sig$ in $\RE^3$. Set $e_+=e_1+e_2+e_3$.
Let the top-dimensional cones be spanned by $\{e_1,e_2,e_+\}$,
$\{e_1,e_2,-e_3\}$, $\{-e_1,-e_2,e_3\}$, $\{-e_1,-e_2,-e_+\}$,
and those obtained by cyclic permutations of the basis.
Define a piecewise linear function $\vph$ on $\Sig$ by
$\vph(e_+)=\vph(-e_i)=1$, $\vph(-e_+)=\vph(e_i)=-1$ ($i=1,2,3$).
Then $\vph$ determines a line bundle such that $0$ is a weight of
both $H^1$ and $H^2$.\\
\noindent
3. It would be interesting to investigate the general conditions 
on the fan under which any weight is of multiplicity $1$
or can appear in only one of the cohomology groups.
One knows that if $M_1$ and $M_2$ are toric manifolds of one of these types,
then so is any $M_1$-fibered toric manifolds over $M_2$.
An interesting class of examples is the Bott-Samelson manifolds
studied from the symplectic point of view in [\ref{GK}].}
\end{rmk}

\sect{Necessity of the \ka condition}

It is well-known that the Hirzebruch-Riemann-Roch theorem or it equivariant
counterpart, the Atiyah-Bott fixed-point theorem, are valid for holomorphic
vector bundles over arbitrary compact complex manifolds.
Therefore it comes as a surprise that the strong equivariant holomorphic
Morse inequalities holds for \ka manifolds only. 
In this section, we show that the strong inequalities (\ref{strong}) are
violated on a suitable non-\ka toric manifold.
We also use the violation of the strong inequalities to show 
the non-existence of invariant \ka structures on certain symplectic manifolds.

Recall that a toric manifold of complex dimension $n$ is \ka if and only if
its fan $\Sig$ admits strictly convex piecewise linear function.
When the complex dimension $n=2$, such a convex function always exists.
Therefore every toric $2$-manifold is K\"ahler.
Notice that in this dimension, the strong inequalities (\ref{strong}) are
equivalent to the weak ones (\ref{weak}).
When $n=3$, there are fans in $\gt\cong\RE^3$ which does not admit any
convex piecewise linear functions.
An example is the fan with $8$ top-dimensional cones $\sig_i$ ($0\le i\le 6$)
and $\sig'$, whose stereographic projection is given by 
Figure~\ref{counter}(a) (see for example [\ref{Da}, \ref{F}]).
Suppose that there were a strictly convex piecewise linear function $\vph$
on $\Sig$.
For the choice of coordinates in Figure~\ref{counter}(a),
convexity on the adjacent cones $\sig_1$ and $\sig_4$ means that
$2\vph(a)+\vph(f)>\hf\vph(c)+\hf\vph(d)$.
Summing over this and the other two inequalities obtained by the 
cyclic permutations of $(a,d)$, $(b,e)$, $(c,f)$, we get 
$3(\vph(a)+\vph(b)+\vph(c))+\vph(d)+\vph(e)+\vph(f)>0$.
On the other hand, convexity on $\sig_4$ and $\sig_2$ implies that
$-2\vph(a)>\vph(b)+\vph(f)$, and hence 
$3(\vph(a)+\vph(b)+\vph(c))+\vph(d)+\vph(e)+\vph(f)<0$, a contradiction.
The corresponding toric variety has an orbifold singularity
because the cone $\sig'$ is not spanned by a $\ZZ$-basis.
(In general, a toric variety has at most orbifold-type singularities if
all the top-dimensional cones in the fan are simplicial.
It is smooth if these cones are spanned by a $\ZZ$-basis.)
The singularity can be avoided by a further triangulation of $\sig'$
into $15$ cones $\sig_i$ ($7\le i\le 21$) by Jurkiewicz [\ref{J}],
as shown in Figure \ref{counter}(b).
This results a smooth toric $3$-manifold with $22$ fixed points.

\begin{prop}\label{CE}
There exists a $T^3$-invariant line bundle over the smooth toric 3-manifold
corresponding to the fan of Figure \ref{counter} such that the strong
equivariant holomorphic Morse inequalities are not satisfied.
\end{prop}
\proof{Let $\{e_1,e_2,e_3\}$ be the standard basis in $\gt\cong\RE^3$,
and $\{e^*_1,e^*_2,e^*_3\}$, its dual basis in $\gt^*$.
We choose a piecewise linear function $\vph$ on the fan 
in Figure \ref{counter} such that
$\vph=-(e^*_1+e^*_2+e^*_3)$ in the cones $\sig_i$, $7\le i\le 21$
(those inside $\sig'$), $\vph=0$ in $\sig_0$, and 
$\vph$ is given by linear interpolations between the cones $\sig'$
and $\sig_0$, for example, $\vph|_{\sig_1}=-3e^*_2+3e^*_3$,
$\vph|_{\sig_4}=-3e^*_3$.
Choose the action chamber $C=\sig_7$ (shown in Figure~\ref{counter}(b)).
It is straightforward to check that\\
1) $\;0\in\Gam^{k,C}$ for $k=0,7$ only and $n^{p_0,C}=0$,
$n^{p_7,C}=3$,\\
2) $\;0\not\in\Gam^{k,-C}$ for any $k$.\\
If the strong inequalities (\ref{strong}) were to hold, then according to 
Proposition \ref{SUPP}, the above claims imply that\\
1) $\;0$ is a weight in $H^0$ and $H^3$ with multiplicity $1$,\\
2) $\;0$ is not a weight in $H^0$ or in $H^3$,\\
respectively. This is a contradiction.}

We actually showed that the strong inequalities (\ref{strong+}) are violated
for the action of any circle subgroup of $T^3$ whose generator is in 
$C=\sig_7$.
The example is consistent with the Atiyah-Bott fixed-point formula (\ref{ab})
since possible contributions to $H^0$ and $H^3$ cancel 
in the equivariant index.

\begin{rmk}
{\em In a forthcoming work [\ref{Wu}], the author shows that the strong
holomorphic Morse inequalities hold for any meromorphic group action on 
a complex manifold with a filterable Bialynicki-Birula decomposition,
which is a filtration of the manifold by closed subvarieties compactible 
with the group action.
This condition is weaker than the \ka assumption made in the current paper. 
The toric manifold in Proposition~\ref{CE} was initially constructed to show 
that Bialynicki-Birula decompositions need not be filterable [\ref{J}].
Therefore the work of [\ref{Wu}] provides a deeper understanding of 
this toric manifold as a counterexample of the strong inequalities on
a general complex manifold.}
\end{rmk}

We now make some remarks on the weak inequalities.

\begin{rmk}\label{dema}
{\em Proposition~\ref{TORIC} on toric manifolds was derived from the weak 
inequalities (\ref{weak}), together with the Atiyah-Bott fixed-point formula
(\ref{ab}).
For a general (possibly non-K\"ahler) toric variety $M$ with a $T$-invariant 
holomorphic line bundle $L_\vph$, Demazure's result [\ref{De}]
(see also [\ref{Da}, \ref{O}, \ref{F}]) states that 
the multiplicity of $\xi\in\lat^*$ in $H^k$ is
equal to the dimension of the local cohomology group $H^k_{Z(\xi)}(\gt)$,
where $Z(\xi)=\set{x\in\gt}{\xi(x)\ge\vph(x)}$.
In fact, Demazure's result implies Proposition \ref{TORIC}.
($\xi\in\lat^*\cap\Gam^0$ if and only if $Z(\xi)=\gt$, in which case
$H^0=\CO$ and $H^k=0$ for $k>0$.
That the multiplicity in $H^0$ is not greater than $1$ also follows from
the existence of a dense open $\TC$-orbit in $M$.)
This suggests that (at least in the completely integrable cases)
the weak holomorphic Morse inequalities could be valid for a larger class
of complex manifolds.
For example, the weak inequalities are not violated for the line bundle
in Proposition~\ref{CE}.
Notice also that the lowest complex dimensions for a toric manifold to be
non-\ka is $3$; this is when 
the weak Morse inequalities become {\em weaker} than the strong ones.
This leaves open the possibility that weak inequalities might hold
when the manifold is complex but not K\"ahler.}
\end{rmk}

\begin{rmk}\label{CONDITION}
{\em More generally, consider a complex manifold $M$ 
with a holomorphic torus $T$-action.
Assume that the fixed-point set $F$ is discrete.
For any $p\in F$, let $\lam^p_1,\cdots,\lam^p_n\in\lat^*\m\{0\}$ be
the weights of $T$ on $T_pM$.
Let $E$ be a holomorphic vector bundle over $M$ with a lifted $T$-action.
Let $e^p_1,\cdots,e^p_r$ form a basis of $E_p$, with weights
$\mu^p_1,\cdots,\mu^p_r\in\lat^*$, respectively.
Near a fixed point $p\in F$, the $T$-action on $M$ has the local model
$\CO^n$ equipped with the action $\e{\ii\tht}\colon(z_1,\cdots,z_n)\mapsto
(\e{\ii\bra\lam^p_1,\tht\ket}z_1,\cdots,\e{\ii\bra\lam^p_n,\tht\ket}z_n)$,
$\tht\in\gt$.
Holomorphic sections of $E$ in a neighborhood of $p$ correspond to
linear combinations of $z^{m_1}_1\cdots z^{m_n}_n\otimes e^p_i$
($m_1,\cdots,m_n\in\NN$, $1\le i\le r$), whose weight under the $T$-action
is $\mu^p_i-\sumk m_k\lam^p_k$.
(The minus sign is explained in [\ref{MW}] in $S^1$-cases.)
Not all such local sections extend to $M$.
So $\supp(H^0)\subset\set{\mu^p_i-\sumk m_k\lam^p_k}{m_1,\cdots,m_n\in\NN,
1\le i\le r}$, counting multiplicities.
This is part of the weak Morse inequalities.
For similar reasons, other weak inequalities in (\ref{weak}) are likely to be
true when the manifold $M$ admits a $T$-invariant complex structure.}
\end{rmk}

Violation of strong holomorphic Morse inequalities can be used as
an obstruction to the existence of $T$-invariant \ka structures.
Recently, Tolman [\ref{T}] constructed a simply-connected six dimensional
symplectic manifold $(M,\om)$ that has a Hamiltonian $T^2$-action with
isolated fixed points but does not admit any $T^2$-invariant \ka structure.
(See [\ref{Wo}] for another perspective.)
Let $\{e_1,e_2\}$ be the standard basis in $\gt\cong\RE^2$ 
and $\{e^*_1,e^*_2\}$, the dual basis in $\gt^*$.
The six fixed points $p_i$ ($1\le i\le 6$) in $M$ correspond to the six
vertices $\mu_i=\mu(p_i)$ of the moment polytope in Figure \ref{tolman}(a)
(reproduced from [\ref{T}, Figure 2]).
Their isotropy weights are the negatives of
$\{e^*_1,e^*_2,e^*_1+e^*_2\}$, $\{e^*_1,e^*_2,-e^*_1-e^*_2\}$,
$\{e^*_1,-e^*_2,e^*_1-e^*_2\}$, $\{-e^*_1,-e^*_2,e^*_1-e^*_2\}$,
$\{-e^*_1,-e^*_1+e^*_2,-2e^*_1+e^*_2\}$,
$\{-e^*_1,-e^*_1+e^*_2,2e^*_1-e^*_2\}$, respectively.
It can be shown that the strong inequalities (\ref{strong}) do not hold
in this example, therefore giving an alternative proof of the non-existence
of $T^2$-invariant \ka structures.
(Most of this comes up in a discussion with S.\ Tolman and C.\ Woodward.)

\begin{prop} {\em (Tolman [\ref{T}])}
There is no $T^2$-invariant \ka structure on Tolman's manifold M.
\end{prop}

\proof{There is a $T^2$-invariant line bundle $L$ such that the weights on
$L_{p_i}$ ($1\le i\le6$) are $\xi_1=0$, $\xi_2=3e^*_1+3e^*_2$, $\xi_3=2e^*_2$,
$\xi_4=3e^*_1+2e^*_2$, $\xi_5=5e^*_1$, $\xi_6=-e^*_1+3e^*_2$
(Figure~\ref{tolman}(b), from [\ref{T}, Figure 3, Case $0<t<s$]).
If there is a $T^2$-invariant \ka form $\om'$, which may be different from
the original symplectic form $\om$, then Figure~\ref{tolman}(b) is
a deformation of the moment polytope associated to $\om'$.
Since all the symplectic quotients of $M$ by the $T^2$-action are $\CO P^1$,
where the complex structure is unique,
the curvature of $L$ is still a $(1,1)$-form on $M$.
Therefore $L$ can be equipped with a holomorphic structure,
and it is invariant under the $T^2$-action.
Let $C$ be the cone spanned by $\{e_1-e_2,-e_2\}$.
Then $e^*_1+2e^*_2\in\Gam\pC$ for $p=p_1,p_5$ only and $n^{p_1,C}=1$,
$n^{p_5,C}=0$ (Figure~\ref{tolman}(c)).
If the strong inequalities (\ref{strong}) were to hold,
by part 1 of Proposition~\ref{SUPP}, $e^*_1+2e^*_2\not\in\supp(H^2)$.
On the other hand, $e^*_1+2e^*_2\in\Gam^{p,-C}$ only for $p=p_3,p_6$
and $n^{p_3,-C}=0$, $n^{p_6,-C}=2$ (Figure~\ref{tolman}(d)).
By part 2 of Proposition~\ref{SUPP}, $e^*_1+2e^*_2\in\supp(H^2)$,
a contradiction.}

In fact, this argument shows further that there is no \ka structure invariant
under any $S^1$-subgroup of $T^2$ whose generator lies in the cone $C$.

\begin{rmk}{\em It is not known whether there is a $T^2$-invariant
complex structure on Tolman's manifold.
If it exists, then the strong inequalities (\ref{strong}) need not be true
even when there is a moment map.
On the other hand since the toric 3-manifold in Proposition~\ref{CE} is not
symplectic, though generalized moment maps in the sense of [\ref{KT}] do 
exist, there is a possibility that the strong inequalities (\ref{strong})
could be true for Hamiltonian torus actions on symplectic manifolds with
invariant complex structures not necessarily {\em calibrated} 
[\ref{Au}, \S II.1.5] by the symplectic forms.
If so, then Tolman's example does not admit any $T^2$-invariant
complex structure.}
\end{rmk}

\newpage

\sect{Non-Abelian equivariant holomorphic Morse inequalities}

In this section we consider a compact \ka manifold $M$ of
complex dimension $n$ with a holomorphic and Hamiltonian action of
a compact connected real Lie group $G$, assuming that the fixed-point
set of a maximal torus $T\subset G$ is discrete.
Let $E$ be a holomorphic vector bundle over $M$ with a lifted holomorphic
$G$-action.
Let $\GC$ be the complexification of $G$ that contains $G$ as 
a maximal compact subgroup.
Then $\GC$ acts holomorphically on $M$ and $E$ [\ref{GS2}].
Choose an (isolated) $T$-fixed-point $p\in F$ and denote the isotropy
group of $p$ in $G$ and $\GC$ by $H=G_p$ and $U=(\GC)_p$, respectively.
(In this way, we have for simplicity dropped the subscript $p$
unless ambiguity occurs.)
Clearly $H\supset T$ and $U$ is a complex subgroup of $\GC$
containing $\HC$ such that $U\cap G=H$.
We need a few group-theoretic results on $\GC$, $U$, $\HC$ satisfying
the above constraints.
These results will determine the behavior of the fixed-point set,
the isotropy weights and the fiber of the vector bundle over the fixed points.
For this purpose, it suffices that $M$ is a complex manifold with
a holomorphic $G$-action; the \ka condition will be used only to
establish the equivariant holomorphic Morse inequalities.

Let $H_0$, $U_0$ be the identity components of $H$ and $U$,
and let $W=N_G(T)/T=N_{\Gc}(\TC)/\TC$, $W_H$, $W_{H_0}$, $W_U$, $W_{U_0}$
be the Weyl groups of $G$, $H$, $H_0$, $U$, $U_0$, respectively.

\begin{lemma}\label{ORBITS}
$W$ acts transitively on the $\TC$-fixed-point set in $\GC/U$
with isotropy group $W_U$.
\end{lemma}

\proof{The fixed-point set is $\set{gU}{g\in \GC, g^{-1}\TC g\subset U}$.
It has a well-defined action of $W$ because $h\in N_{\Gc}(\TC)$
acts on the fixed points by $gU\mapsto hgU$ and for any $t\in \TC$,
$htgU=hg(g^{-1}tg)U=hgU$.
To show that the action is transitive, choose any fixed point $gU$.
Since $g^{-1}\TC g$ and $\TC$ are two maximal tori in $U_0$,
there is an element $u\in U_0$ such that $u^{-1}g^{-1}\TC gu=\TC$.
So $gu\in N_{\Gc}(\TC)$ and $gU=(gu)U$ can be obtained from $U$ by
an element of $W$.
Finally, $h\in N_{\Gc}(\TC)$ fixes $U$ (i.e., $hU=U$)
if and only if $h\in U\cap N_{\Gc}(\TC)=N_U(\TC)$.
So the isotropy group in $W$ is $N_U(\TC)/\TC=W_U$.}

\begin{prop}
There is an action of $W$ on the $T$-fixed-point set $F$ in $M$.
Each $W$-orbit in $F$ is the intersection of $F$ with a $G$- or
$\GC$-orbit in $M$, which do not depend on the choice of $T$.
\end{prop}

\proof{Lemma~\ref{ORBITS} implies that for $p\in F$, the $T$-fixed points
in $\GC\cdot p\cong\GC/U$ form a single $W$-orbit.
For a different maximal torus $gTg^{-1}$ ($g\in G$), the new fixed-point set
$gF$ is contained in the same $G$-orbits.}

\begin{ex}
{\em We consider the diagonal action of $G=SO(3)$ on $M=S^2\times S^2$,
where $G$ acts on $S^2$ by standard rotations.
The maximal torus in $G$ is $T=S^1$.
Let $n$, $s$ be the poles in $S^2$ which are fixed by $T$.
Then the $T$-fixed-point set in $M$ is $F=\{(n,n),(s,s),(n,s),(s,n)\}$.
Though the isotropy groups $G_p=T$ for all $p\in F$, the situation is different
if we consider the action of the complexification $\GC=SL(2,\CO)/\ZZ_2$.
It is easy to see that the $(\GC)_{(n,n)}=B$, a Borel subgroup in $\GC$,
whereas $(\GC)_{(n,s)}=\TC$.
Consequently, the orbit $\GC\cdot(n,n)=\set{(x,x)}{x\in S^2}\cong\GC/B$ is
compact and contains another fixed point $(s,s)$ related to $(n,n)$ by
$W=\ZZ_2$, whereas $\GC\cdot(n,s)=\set{(x,y)}{x\ne y\in S^2}\cong\GC/\TC$
is non-compact and contains $(s,n)$, related to $(n,s)$ also by $\ZZ_2$.
In this case, $M$ is the union of two $\GC$-orbits.}
\end{ex}

Let $\Del$, $\Del_{H_0}$, $\Del_{U_0}$ be the set of roots of
the pairs $(\GC,\TC)$, $(\HC_0,\TC)$, $(U_0,\TC)$, respectively.
Choose a set of positive roots $\Del^+\subset\Del$ and let
$\Del^-=-\Del^+$, $\Del^\pm_{H_0}=\Del^\pm\cap\Del_{H_0}$.
The length of $w\in W$ is $l(w)=|w\Del^+\cap\Del^-|$.
Let $l_{H_0}(w)=|w(\Dp_{H_0})\cap\Del^-|$ be the length of $w$
{\em relative} to the subgroup $H_0$.
Notice that $l_T(w)=l(w)$ ($w\in W$).

\begin{lemma}\label{ISOTROPY}
There is a $\TC$-fixed-point in $\GC/U$ at which the set of weights of the
isotropy representation contains $\Dp_{H_0}$.
\end{lemma}

\proof{Since the sets of isotropy weights at two fixed-points
$p$ and $wp$ ($w\in W$) are related by the action of $w$ and
since any two choices of $\Del^+$ are conjugate under a $w\in W$,
it suffices to show that the lemma holds for any given $p$ under
a particular choice of $\Del^+$ depending on $p$.
Because $U$ is a complex subgroup of $G$, $\Del_{U_0}$ is a {\em closed}
subset of $\Del$, i.e., if $\al,\beta\in\Del_{U_0}$ and $\al+\beta\in\Del$,
then $\al+\beta\in\Del_{U_0}$.
Furthermore, since $U\cap G=H$, $\Del_{U_0}\cap(-\Del_{U_0})=\Del_{H_0}$.
It is easy to see that $\Del'=\Del_{U_0}\m\Del_{H_0}$ is also closed
(if $\al,\beta\in\Del'$ but $\al+\beta\in\Del_{H_0}$,
then $-\beta=(-\al-\beta)+\al\in\Del_{U_0}$, a contradiction)
and satisfies $\Del'\cap(-\Del')=\emptyset$.
So $\Del'$ is contained in some choice of $\Del^-$ [\ref{Bou}].
Consequently, there is a parabolic subgroup $P\supset U_0$ and $P\cap G=H_0$.
Therefore the set of weights of isotropy representation at $p$ is
$\Del\m\Del_{U_0}\supset\Del\m\Del_P=\Dp_{H_0}$.}
\vspace{2ex}

For each $W$-orbit $S\in F/W$, choose $p\in S$ and let $H=G_p$ be the 
isotropy group in $G$.
Denote $\Del_S=\Del_{H_0}$, $l_S(w)=l_{H_0}(w)$ and
$\det_Sw=(-1)^{l_S(w)}\det w$ for $w\in W$;
they do not depend on the choice of $p$ in $S$.

\begin{prop}\label{POSITIVE}
In each orbit $S\in F/W$ we can choose a representative, denoted by $p_S$,
such that the weights of the isotropy representation at $p_S$
contains $\Dp_S$.
\end{prop}

\proof{This is a direct consequence of Lemma~\ref{ISOTROPY}.}

\begin{rmk}
{\em In [\ref{GP}] (see also [\ref{PW}]), under a regularity assumption
that implies $H=T$, it was shown that the set of isotropy weights contains 
$\Del^+$ up to signs that depend on the roots.
Proposition~\ref{POSITIVE} is a refinement of this result when $M$ is complex
and when the $G$-action preserves the complex structure.
Let $\lSk$ ($k=1,\cdots,n-|\Dp_S|$) be the other weights at $p_S$.
Since $T_{p_S}M$ is also a representation of $U$, the set of weights
$(\Dp_S)\cup\{\lSk\}$ is $W_U$-invariant.
At another fixed point $wp_S\in S$,
the corresponding set of isotropy weights is $w((\Dp_S)\cup\{\lSk\})$,
which depends only on the coset $\wb\in W/W_U$.
By $W$-invariance, any action chamber can be transformed into one that
intersects the positive Weyl chamber $\gt^+=\set{x\in\gt}{\bra\Del^+,x\ket>0}$
in an open cone; let $C$ be one of such.
At $p_S$, the polarizing index $n^C_{p_S}$ is equal to $n^C(\{\lSk\})$,
that of the set $\{\lSk\}$.
At $wp_S$, $n^C_{wp_S}=l_S(w)+n^C(\{\wl\})$.}
\end{rmk}

The following lemma and example are communicated to the author by D.\ Vogan.

\begin{lemma}\label{COMPONENTS}
$N(\Del_{U_0})=\set{w\in W}{w\Del_{U_0}=\Del_{U_0}}$ contains $W_{H_0}$ as
a normal subgroup.
Moreover, $U/U_0\cong W_U/W_{H_0}$ is a subgroup of $N(\Del_{U_0})/W_{H_0}$.
\end{lemma}

\proof{Since $\HC_0$ is a subgroup of $U$, $W_{H_0}$ and $W_U$ preserve 
the set $\Del_{U_0}$, hence $W_{H_0}<W_U<N(\Del_{U_0})$.
For any $\al\in\Del_{H_0}$, $w\in N(\Del_{U_0})$,
since $wr_\al w^{-1}=r_{w\al}$, and 
$w\al\in w\Del_{U_0}\cap(-w\Del_{U_0})=\Del_{H_0}$,
$W_{H_0}$ is a normal subgroup of $N(\Del_{U_0})$.
Finally, we prove $U/U_0\cong W_U/W_{H_0}$.
Define a homomorphism $N_U(\TC)\to U/U_0$ by $u\mapsto uU_0$.
Its kernel is $N_{U_0}(\TC)$; we show that it is onto.
For any $u\in U$, $u^{-1}\TC u$ is a maximal torus in $U_0$.
So there is $u_0\in U_0$ such that $u^{-1}_0u^{-1}\TC uu_0=\TC$.
This implies that $uu_0\in N_U(\TC)$ and $uu_0\mapsto uU_0$.
By homomorphism theorems, $U/U_0\cong N_U(\TC)/N_{U_0}(\TC)\cong W_U/W_{U_0}$.
Finally, let $P$ be a parabolic subgroup that contains $U_0$ and such that 
$P\cap G=H_0$, then $W_{H_0}<W_{U_0}<W_P$.
Since $W_P=W_{H_0}$, we conclude that $W_{U_0}=W_{H_0}$.}

\begin{ex}
{\em The isotropy group $U$ can be disconnected.
Let the action of $G=SO(3)$ on $M=\CO P^2$ be the projectivization of the
adjoint representation on $\mbox{\frak{so}}\,(3)^\co\cong\CO^3$.
If the root space decomposition is
$\mbox{\frak{so}}\,(3)^\co=\CO e_0\oplus\CO e_+\oplus\CO e_-$,
where $\CO e_0=\mbox{Lie}(\TC)$ and $[e_0,e_\pm]=\pm e_\pm$,
then the $T$-fixed-point set $F=\{[e_0],[e_+],[e_-]\}$.
The isotropy groups of $[e_\pm]$ in $\GC$ are Borel subgroups
while that of $[e_0]$ is the full normalizer $N_{\Gc}(\TC)$ of $\TC$,
which is disconnected.}
\end{ex}

Since $G$ is compact, $\gt$ is equipped with an invariant inner product;
this induces one on $\gt^*$, which will be denoted by $(\cdot,\cdot)$.
Let $\dom=\set{\lam\in\lat^*}{(\lam,\Del^+)\ge0}$ and
$\dom_{H_0}=\set{\lam\in\lat^*}{(\lam,\Del^+_{H_0})\ge0}$ be the sets of
dominant weights with respect to $G$ and $H_0$, respectively.
We denote the irreducible representation of $G$ ($H_0$, respectively) of
the highest weight $\lam\in\dom$ ($\lam\in\dom_{H_0}$, respectively) by
$R^G_\lam$ ($R^{H_0}_\lam$, respectively).
Since $W_U$ preserves $\Del_{U_0}$, hence
$\Del_{U_0}\cap(-\Del_{U_0})=\Del_{H_0}$, the action of $W_U$
permutes the Weyl chambers of the pair $(\HC_0,\TC)$.
Let $\dom_U\subset\dom_{H_0}$ be a fundamental region of the group 
$W_U$ in $\lat^*$.

\begin{lemma}\label{REPR}
Let $V$ be a finite dimensional representation of $U$, then
	\be\label{char}
\ch(V)=\sum_{\Lam\in\dom_U}m_\Lam\sum_{w\in W_U}w\,\frac{\e{\Lam}}{\pra{H_0}}.
	\ee
Here $m_\Lam\in\QQ$; if $U$ is connected, then $m_\Lam\in\NN$
is the multiplicity of $R^H_\Lam$ in $V$.
\end{lemma}

\proof{As a representation of $H_0$, 
$V=\bigoplus_\lam R^{H_0}_\lam$; the direct sum is over some
$\lam\in\dom_{H_0}$ with possible multiplicities.
Since $V$ is a representation of $U$, $W_U/W_{H_0}$ acts on the set
$\{R^{H_0}_\lam\}$.
For a given $W_U/W_{H_0}$-orbit through $R^{H_0}_\Lam$,
let $W_\Lam$ be the isotropy group.
The contribution of this orbit to the character is
	\be
\inv{|W_\Lam|}\sum_{\wb\in W_U/W_{H_0}}\wb\ch(R^{H_0}_\Lam)
=\inv{|W_\Lam|}\sum_{w\in W_U}w\,\frac{\e{\Lam}}{\pra{H_0}}.
	\ee
Equation (\ref{char}) follows readily.
If $U$ is connected, then $W_U/W_{H_0}=\{1\}$, hence all $m_\Lam\in\NN$.}
\vspace{2ex}

At $p_S$ ($S\in F/W$) whose isotropy group in $\GC$ is $U$,
set $\dom_S=\dom_U$.
According to the above lemma,
	\be\label{subgroup}
\ch(E_{p_S})=\sum_{\Lam\in\dom_S}m^S_\Lam
\sum_{w\in W_U}w\frac{\e{\Lam}}{\pra{S}},
	\ee
where $m^S_\Lam\in\QQ$; if the isotropy group $U$ is connected,
then $m^S_\Lam\in\NN$ is the multiplicity of $R^H_\Lam$ in $E_{p_S}$.

\begin{thm}\label{NONAB}
Let $G$ be a compact connected Lie group acting Hamiltonianly on
a compact \ka manifold $M$ preserving the complex structure.
Assume that the fixed-point set $F$ of the maximal torus $T$ is discrete.
If $E$ is a holomorphic  vector bundle over $M$ where the $G$-action 
lifts holomorphically,
then $G$ acts on the cohomology groups $H^k(M,\OO(E))$ ($0\le k\le n$);
let $H^k(M,\OO(E))=\bigoplus_{\Lam\in\dom}m\kL R^G_\Lam$ ($m\kL\in\NN$)
be the decompositions into irreducible representations $R^G_\Lam$ of $G$.
Choose an action chamber $C$ that intersects $\gt^+$ in an open cone.
Then under the above notations, we have the following non-Abelian
equivariant holomorphic Morse inequalities
	\bea\label{nonab}
\vc\sum_{S\in F/W}\sumw\det\!_Sw\;t^{l_S(w)+\nwl}
   \sum_{\Lam\in\dom_S}m\SL\,\e{\wL}\!\prwl				\nno
\eq\sumk t^k\sumw\det w\,\sum_{\Lam\in\dom}m\kL\,\e{\wL}+(1+t)Q^C(t)\pra{}
	\eea
for some $Q^C(t)\ge 0$.
\end{thm}

\proof{We apply (\ref{strong}) to the maximal torus $T$ of $G$.
The contribution of one orbit $S\in W/F$ is a sum over cosets $\wb\in W/W_U$.
This sum can be combined with that over $W_U$ in (\ref{subgroup}).
Since the set $w((\Dp_S)\cup\{\lam^S_k\})$ depends only on the coset
$\wb\in W/W_U$ and since $\ch(E_{\wb p_S})=\wb\ch(E_{p_S})$, the result is
	\bea
\vc\sumw t^{n^C_{wp_S}}\prod_{w\al\in w(\Dp_S)\cap\Del^+}\inv{1-\e{-w\al}}
   \prod_{w\al\in w(\Dp_S)\cap\Del^-}\frac{\e{w\al}}{1-\e{w\al}}
   \sum_{\Lam\in\dom_S}m\SL\,w\left(\frac{\e{\Lam}}{\pra{S}}\right)
   L^C_{w,S}								\nno
\eq\sumw(-1)^{l(w)-l_S(w)}t^{l_S(w)+\nwl}\,w
   \left(\prod_{\al\in\Dp_S}\inv{1-\e{-\al}}\right)
   \sum_{\Lam\in\dom_S}m\SL w\left(\frac{\e{\Lam}}{\pra{S}}\right)L^C_{w,S}\nno
\eq\sumw\det\!_Sw\;t^{l_S(w)+\nwl}
   \frac{\det w\sum_{\Lam\in\dom_p}m\SL\,\e{\wL}}{\pra{}}\,L^C_{w,S},
	\eea
where 
	\be
L^C_{w,S}=\prwl
	\ee
is a factor from weights not in $\Dp_S$.
Using Weyl's character formula, the contribution to (\ref{strong})
from the cohomology group is
	\be
\sumk t^k\ch(H^k)=
\sumk t^k\sum_{\Lam\in\dom}m\kL\frac{\sumw\det w\,\e{\wL}}{\pra{}}.
	\ee
The result follows after multiplying both sides of (\ref{strong}) by $\pra{}$.}
\vspace{2ex}

\begin{rmk}
{\em Setting $t=-1$ in (\ref{nonab}), we obtain a fixed-point formula
	\be\label{nonabfix}
\sum_{S\in F/W}\sumw\det w\sum_{\Lam\in\dom_S}m\SL\,\e{\wL}\!
\prod_{k=1}^{n-|\Dp_S|}\inv{1-\e{-w\lam^S_k}}	
=\sumk(-1)^k\sumw\det w\,\sum_{\Lam\in\dom}m\kL\,\e{\wL}
	\ee
for non-Abelian group actions.}
\end{rmk}

\vskip4ex

\noindent {\bf Acknowledgement.} 
The author would like to thank M.\ Bozicevic, R.\ Dabrowski, M.\ Grossberg,
E.\ Prato, S.\ Tolman, D.\ Vogan, C.\ Woodward for illuminating discussions.
He also thanks M.\ Vergne for kindly pointing out an error in an earlier
version of this paper and the referees for helpful suggestions.
This work is partially supported by NSF grants DMS-93-05578 at Columbia
University, DMS-90-22140 at MSRI and by ICTP.

\bigskip

        \newcommand{\athr}[2]{{#1}.\ {#2}}
        \newcommand{\au}[2]{\athr{{#1}}{{#2}},}
        \newcommand{\an}[2]{\athr{{#1}}{{#2}} and}
%       \au{first name}{last name}
        \newcommand{\jr}[6]{{#1}, {\it {#2}} {#3}\ ({#4}) {#5}-{#6}}
%       \jr{title}{journal}{volume}{year}{first-page}{last-page}
        \newcommand{\pr}[3]{{#1}, {#2} ({#3})}
%       \pr{title}{preprint number}{year}
        \newcommand{\bk}[4]{{\it {#1}}, {#2}, ({#3}, {#4})}
%       \bk{title}{publisher}{place}{year}
        \newcommand{\cf}[8]{{\it {#1}}, {#2}, {#5},
                 {#6}, ({#7}, {#8}), pp.\ {#3}-{#4}}
%\cf{title}{conf-title}{f-page}{l-page}{editor}{publisher}{place}{year}
%      1        2          3      4       5         6        7     8
        \vspace{5ex}
        \begin{flushleft}
{\bf References}
        \end{flushleft}
{\small
        \begin{enumerate}
        
        \item\label{A}
        \au{M.\ F}{Atiyah}
        \jr{Convexity and commuting Hamiltonians}
        {Bull. London Math. Soc.}{14}{1982}{1}{15}

        \item\label{AB}
        \an{M.\ F}{Atiyah} \au{R}{Bott}
        \jr{A Lefschetz fixed point formula for elliptic complexes, Part I}
        {Ann. Math.}{86}{1967}{374}{407};
        \jr{Part II}{Ann. Math.}{87}{1968}{451}{491}

	\item\label{Au}
	\au{M}{Audin}
	\bk{The topology of torus actions on symplectic manifolds,
	{\rm Prog. in Math. 93}}
	{Birkh\"auser Verlag}{Basel, Boston, Berlin}{1991}

%	\item\label{B}
%	\au{R}{Bott}
%	\jr{Homogeneous vector bundles}
%	{Ann. Math.}{66}{1957}{203}{248}

	\item\label{Bou}
	\au{N}{Bourbaki}
	\bk{Groupes et alg\`ebres de Lie, \'El\'ements de Math\'ematique XXXIV}
	{Hermann}{Paris}{1968}, ch.\ 6, p.\ 163

	\item\label{Da}
	\au{V.\ I}{Danilov}
	\jr{The geometry of toric varieties}
	{Russian Math. Surveys}{33}{1978}{97}{154}

	\item\label{Del}
	\au{T}{Delzant}
	\jr{Hamiltonien p\'eriodiques et image convexe de l'application
	moment}{Bull. Soc. Math. France}{116}{1988}{315}{339}

	\item\label{De}
	\au{M}{Demazure}
	\jr{Sous-groupes alg\'ebriques de rang maximum de groupes de Cremona}
	{Ann. Sci. \'Ecole Norm. Sup. (4)}{3}{1970}{507}{588}

	\item\label{Fr}
	\au{T}{Frankel}
	\jr{Fixed points and torsion on \ka manifolds}
	{Ann. Math.}{70}{1959}{1}{8}

	\item\label{F}
	\au{W}{Fulton}
	\bk{Introduction to toric varieties}
	{Princeton Univ. Press}{Princeton, NJ}{1993}

%	\item\label{Gr}
%	\au{Ph.\ A}{Griffiths}
%	\jr{Some geometric and analytic properties of homogeneous complex
%	manifolds, Parts I, II}{Acta Math.}{110}{1963}{115}{208}

	\item\label{GK}
	\an{M}{Grossberg} \au{Y}{Karshon}
	\jr{Bott towers, complete integrability, and the extended character
	of representations}{Duke Math. J.}{76}{1994}{23}{58}

	\item\label{G}
	\au{V}{Guillemin}
	\jr{Kaehler structures on toric varieties}
	{J. Diff. Geom.}{40}{1994}{285}{309}

	\item\label{GLS}
	\au{V}{Guillemin} \an{E}{Lerman} \au{S}{Sternberg}
	\jr{On the Kostant multiplicity formula}
	{J. Geom. Phys.}{5}{1988}{721}{750}

        \item\label{GP}
        \an{V}{Guillemin} \au{E}{Prato}
        \jr{Heckman, Kostant, and Steinberg formulas for symplectic manifolds}
        {Adv. Math.}{82}{1990}{160}{179}

	\item\label{GS}
        \an{V}{Guillemin} \au{S}{Sternberg}
        \jr{Convexity properties of the moment mapping}
        {Invent. Math.}{67}{1982}{491}{513}

	\item\label{GS2}
        \an{V}{Guillemin} \au{S}{Sternberg}
	\jr{Geometric quantization and multiplicities of group representations}
	{Invent. Math.}{67}{1982}{515}{538}

	\item\label{Ho}
	\au{L}{H\"ormander}
        \bk{The analysis of linear partial differential operators I, 2nd ed.}
        {Springer-Verlag}{Berlin, Heidelberg, New York, Tokyo}{1990}{, \S 7.4}

	\item\label{J}
	\au{J}{Jurkiewicz}
	\jr{An example of algebraic torus action which determines the
	nonfiltrable decomposition}{Bull. l'Acad. Polonaise Sci.}
	{25}{1977}{1089}{1092}

	\item\label{KT}
	\an{Y}{Karshon} \au{S}{Tolman}
	\jr{The moment map and line bundles over presymplectic manifolds}
	{J. Diff. Geom.}{38}{1993}{465}{484}

	\item\label{MW}
	\an{V}{Mathai} \au{S}{Wu}
	\jr{Equivariant holomorphic Morse inequalities I: a heat kernel proof}
	{J. Diff. Geom.}{46}{1997}{78}{98}

	\item\label{M}
	\au{Y}{Matsushima}
	\jr{Sur les espaces homog\`enes K\"ahl\'eriens d'un groupe de Lie 
	r\'eductif}{Nagoya Math. J.}{11}{1957}{53}{60}

	\item\label{O}
	\au{T}{Oda}
	\bk{Convex bodies and algebraic geometry}
	{Springer-Verlag}{Berlin, Heidelberg, New York}{1988};
	\cf{Geometry of toric varieties}
	{Proc.\ of the Hyderabad Conference on Algebraic Groups}{407}{440}
	{eds. S.\ Ramanan et.\ al.}{Manoj Prakashan}{Madras}{1991}
	
	\item\label{PW}
	\an{E}{Prato} \au{S}{Wu}
	\jr{Duistermaat-Heckman measures in a non-compact setting}
	{Comp. Math.}{94}{1994}{113}{128}

	\item\label{T}
	\au{S}{Tolman}
	\jr{Examples of non-K\"ahler Hamiltonian torus actions}
	{Invent.\ Math.}{131}{1998}{299}{310}
	
	\item\label{Wa}
	\au{G}{Warner}
	\bk{Harmonic analysis on semi-simple Lie groups I,
	{\rm Grund.\ Math.\ Wiss.\ 188}}
	{Springer-Verlag}{Berlin, Heidelberg, New York}{1992}, ch.\ 1, p.\ 17

	\item\label{W} 
	\au{E}{Witten} 
	\cf{Holomorphic Morse inequalities}
	{Algebraic and differential topology, Teubner-Texte Math., 70}
	{318}{333}{ed.\ G.\ Rassias}{Teubner}{Leipzig}{1984}

	\item\label{Wo}
	\au{C}{Woodward}
	\jr{Multiplicity-free Hamiltonian actions need not be K\"ahler}
	{Invent.\ Math.}{131}{1998}{311}{319}

	\item\label{Wu}
	\au{S}{Wu}
	On the instanton complex of holomorphic Morse theory, in preparation

        \end{enumerate}}
	\newpage

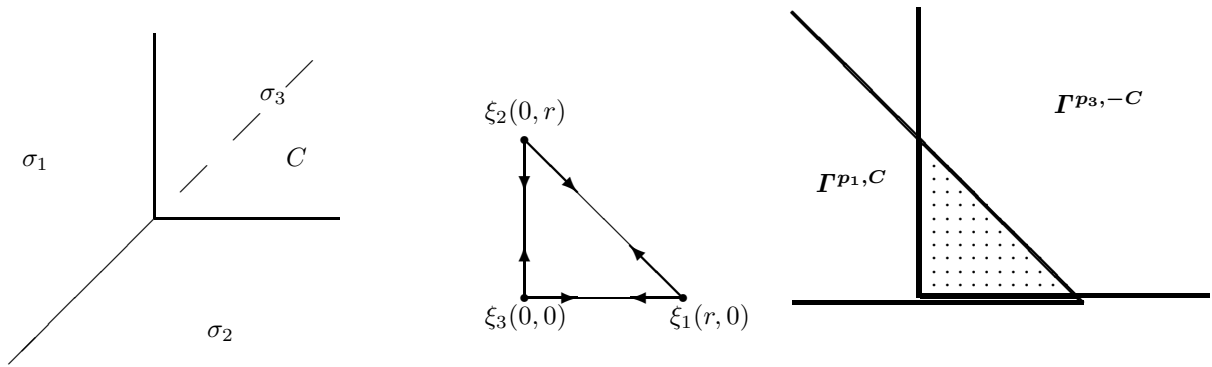
\newcommand{\fig}[2]{	\begin{picture}(20,20)
{#1}
\put(-7.5,-12.5){\parbox{150pt}{{\small {#2}}}}
	\end{picture}}
\newcommand{\lfig}[2]{	\begin{picture}(20,20)
{#1}
\put(-7,-15){\parbox{120pt}{{\small {#2}}}}
	\end{picture}}

\setcounter{sect}{2}	\setcounter{figure}{0}

\begin{figure}[bottom]
	\begin{picture}(450,525)
\unitlength 10pt
\put(12.5,42.5){\fig{\put(0,0){\circle*{.25}}
	\thicklines
\put(0,0){\vector(1,0){5}}
\put(0,0){\vector(0,1){5}}
\put(0,0){\vector(-1,1){5}}
\put(5.5,-0.5){$e^*_1$}
\put(0,5.5){$e^*_2$}
\put(-6,5.5){$e^*_2-e^*_1$}
\put(6,-7){$\gt^*$}}
{(a) An example of isotropy weights at a fixed point $p\in F$ when $\dim T=2$,
$\dim_\co M=3$.}}

\put(32.5,42.5){\fig{\put(-7,-7){\line(1,1){14}}
\put(-7.5,0){\line(1,0){15}}	\put(0,-7.5){\line(0,1){15}}
\put(5,2){$C$}	\put(-6,-2.5){$-C$}	\put(6,-7){$\gt$}}
{(b) The hyperplanes cut $\gt$ into 6 chambers.
(Weights at other fixed points may cut them further.)
A chamber $C$ is chosen.}}

\put(12.5,15){\fig{\put(0,0){\vector(1,0){5}}
\put(0,0){\vector(0,1){5}}		\put(0,0){\vector(-1,1){5}}
\put(0,-1){$\xi$}		\put(-6,-3.5){$\bfgam{p,C}$}
\multiput(-7,6)(-1,0){1}{\circle*{.1}}	\multiput(-6,5)(-1,0){2}{\circle*{.1}}
\multiput(-5,4)(-1,0){3}{\circle*{.1}}	\multiput(-4,3)(-1,0){4}{\circle*{.1}}
\multiput(-3,2)(-1,0){5}{\circle*{.1}}	\multiput(-2,1)(-1,0){6}{\circle*{.1}}
\multiput(-1,0)(-1,0){7}{\circle*{.1}}	\multiput(-1,-1)(-1,0){7}{\circle*{.1}}
\multiput(-1,-2)(-1,0){7}{\circle*{.1}}	\multiput(-1,-3)(-1,0){3}{\circle*{.1}}
\put(-7,-3){\circle*{.1}}		\multiput(-1,-4)(-1,0){7}{\circle*{.1}}
\multiput(-1,-5)(-1,0){7}{\circle*{.1}}	\multiput(-1,-6)(-1,0){7}{\circle*{.1}}
\multiput(-1,-7)(-1,0){7}{\circle*{.1}}	\multiput(-1,-8)(-1,0){7}{\circle*{.1}}
\thicklines	\linethickness{1.5pt}
\put(-0.28,0){\line(-1,1){7}}	\put(-0.33,-.05){\line(-1,1){7}}
\multiput(-0.3,0)(0,-2.3){4}{\line(0,-1){1.15}}}
{(c) The cone $\Gam^{p,C}$ when $\supp E_p=\{\xi\}$.
Dotted line stands for boundaries not included. 
The polarizing index $n^{p,C}=1$.}}

\put(32.5,15){\fig{\put(0,0){\vector(1,0){5}}
\put(0,0){\vector(0,1){5}}		\put(0,0){\vector(-1,1){5}}
\put(-.5,-1){$\xi$}		\put(6,4){$\bfgam{p,-C}$}
\multiput(1,8.5)(1,0){9}{\circle*{.1}}	\multiput(1,7.5)(1,0){9}{\circle*{.1}}
\multiput(1,6.5)(1,0){9}{\circle*{.1}}	\multiput(1,5.5)(1,0){9}{\circle*{.1}}
\multiput(1,4.5)(1,0){5}{\circle*{.1}}	\put(9,4.5){\circle*{.1}}
\multiput(1,3.5)(1,0){9}{\circle*{.1}}	\multiput(1,2.5)(1,0){9}{\circle*{.1}}
\multiput(1,1.5)(1,0){9}{\circle*{.1}}	\multiput(1,.5)(1,0){9}{\circle*{.1}}
\multiput(2,-.5)(1,0){8}{\circle*{.1}}	\multiput(3,-1.5)(1,0){7}{\circle*{.1}}
\multiput(4,-2.5)(1,0){6}{\circle*{.1}}	\multiput(5,-3.5)(1,0){5}{\circle*{.1}}
\multiput(6,-4.5)(1,0){4}{\circle*{.1}}	\multiput(7,-5.5)(1,0){3}{\circle*{.1}}
\multiput(8,-6.5)(1,0){2}{\circle*{.1}}
\thicklines\linethickness{1.5pt}
\multiput(0.2,0)(0,2.6){4}{\line(0,1){1.3}}
\multiput(0.18,0)(2,-2){4}{\line(1,-1){1}}
\multiput(0.23,.05)(2,-2){4}{\line(1,-1){1}}}
{(d) The cone $\Gam^{p,-C}$ when $\supp E_p=\{\xi\}$. 
Dotted line stands for boundaries not included.
The polarizing index $n^{p,-C}=2$.}}
	\end{picture}
\caption{\label{definition} Illustration of isotropy weights, action
chambers and the regions $\Gam^{p,C}$.}
	\end{figure}

\setcounter{sect}{3}	\setcounter{figure}{0}

\newcommand{\projective}{\put(-3,-3){\line(1,0){6}}
\put(-3,-3){\line(0,1){6}}	\put(-3,3){\line(1,-1){6}}}
\begin{figure}[bottom]
	\begin{picture}(450,570)
\unitlength 10pt
\put(7.5,46){\lfig{\multiput(-1,1)(2,2){3}{\line(1,1){1}}
\put(-2,0){\line(1,0){7}}\put(-2,0){\line(0,1){7}}\put(-2,0){\line(-1,-1){5.5}}
\put(2,4.5){$\sig_3$}	\put(-7,2){$\sig_1$}	\put(0,-4.5){$\sig_2$}
\put(3,2){$C$}}
{(a) The fan of $\CO P^2$ with three top dimensional cones $\sig_i$ 
($1\le i\le 3$). An action chamber $C=\sig_3\cap(-\sig_1)$ is chosen.}}

\put(22.5,46){\lfig{\projective		\put(-4.5,-4){$\xi_3(0,0)$}	
\put(2.5,-4){$\xi_1(r,0)$}		\put(-4.5,3.75){$\xi_2(0,r)$}
\put(-3,-3){\circle*{.25}}\put(3,-3){\circle*{.25}}\put(-3,3){\circle*{.25}}
\thicklines
\put(-3,-3){\vector(1,0){2}}		\put(-3,-3){\vector(0,1){2}}
\put(3,-3){\vector(-1,0){2}}		\put(3,-3){\vector(-1,1){2}}
\put(-3,3){\vector(0,-1){2}}		\put(-3,3){\vector(1,-1){2}}}
{(b) Position of weights $\xi_i$ in the fiber over $p_i\in F$.
The two arrows at $\xi_i$ are the negative of isotropy weights at $p_i$.}}

\put(37.5,46){\lfig{\projective
\multiput(-2.5,-2.5)(0,.5){11}{\circle*{.1}}
\multiput(-2,-2.5)(0,.5){10}{\circle*{.1}}
\multiput(-1.5,-2.5)(0,.5){9}{\circle*{.1}}
\multiput(-1,-2.5)(0,.5){8}{\circle*{.1}}
\multiput(-.5,-2.5)(0,.5){7}{\circle*{.1}}
\multiput(0,-2.5)(0,.5){6}{\circle*{.1}}
\multiput(.5,-2.5)(0,.5){5}{\circle*{.1}}
\multiput(1,-2.5)(0,.5){4}{\circle*{.1}}
\multiput(1.5,-2.5)(0,.5){3}{\circle*{.1}}
\multiput(2,-2.5)(0,.5){2}{\circle*{.1}}
\put(2.5,-2.5){\circle*{.1}}
\put(-7,1){$\bfgam{p_1,C}$}		\put(2,4){$\bfgam{p_3,-C}$}
\thicklines	\linethickness{1.5pt}
\put(3.15,-3.17){\line(-1,0){11}}	
\put(3.15,-3.15){\line(-1,1){11}}	\put(3.1,-3.2){\line(-1,1){11}}
\put(-3.1,-2.9){\line(1,0){11}}		\put(-3.1,-3){\line(0,1){11}}}
{(c) The two cones $\Gam^{p_1,C}$ and $\Gam^{p_3,-C}$, with 
$n^{p_1,C}=n^{p_3,-C}=0$, intersect at the dotted region.}}

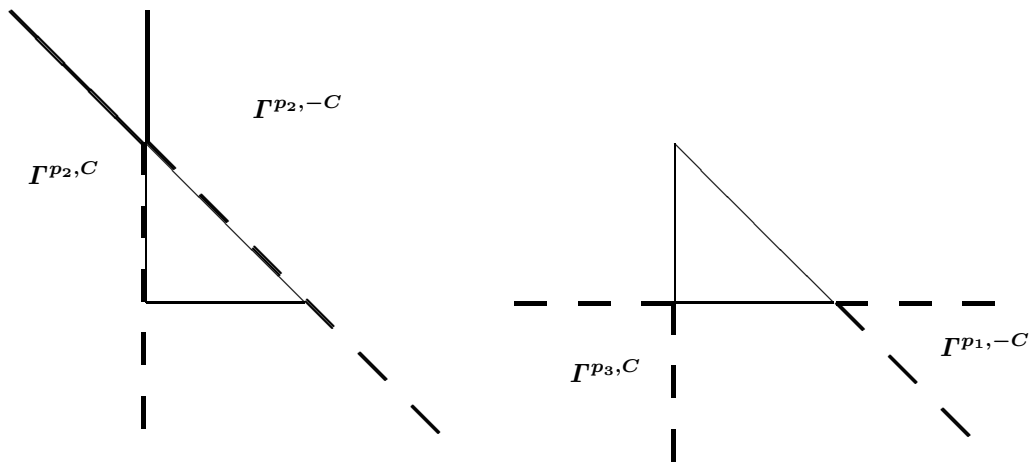
\put(15,18){\lfig{\projective
\put(-7.5,1.5){$\bfgam{p_2,C}$}		\put(1,4){$\bfgam{p_2,-C}$}
\thicklines	\linethickness{1.5pt}
\put(-2.95,3.07){\line(0,1){5}}
\multiput(-2.95,3.05)(2,-2){6}{\line(1,-1){1}}
\multiput(-2.9,3.1)(2,-2){6}{\line(1,-1){1}}
\put(-3.1,3.05){\line(-1,1){5}}	\put(-3.15,3){\line(-1,1){5}}
\multiput(-3.1,3.05)(0,-2.4){5}{\line(0,-1){1.2}}}
{(d) The two cones $\Gam^{p_2,C}$ and $\Gam^{p_2,-C}$, with 
$n^{p_2,C}=n^{p_2,-C}=1$, have empty intersection.}}

\put(35,18){\lfig{\projective
\put(-7,-6){$\bfgam{p_3,C}$}		\put(7,-5){$\bfgam{p_1,-C}$}
\thicklines	\linethickness{1.5pt}
\multiput(-3.05,-3)(0,-2.4){3}{\line(0,-1){1.2}}
\multiput(-3.05,-3.05)(-2.4,0){3}{\line(-1,0){1.2}}
\multiput(3.1,-3.05)(2.4,0){3}{\line(1,0){1.2}}
\multiput(3.1,-3.05)(2,-2){3}{\line(1,-1){1}}
\multiput(3.15,-3)(2,-2){3}{\line(1,-1){1}}}
{(e) The two cones $\Gam^{p_3,C}$ and $\Gam^{p_1,-C}$, with 
$n^{p_2,C}=n^{p_2,-C}=2$, have empty intersection.}}
	\end{picture}
\caption{\label{projective}The projective space $\CO P^2$ with an
invariant line bundle of first Chern number $r\ge 0$}
	\end{figure}

\newcommand{\hirzebruch}{\put(-4,-2){\line(1,0){8}}
\put(-4,-2){\line(0,1){10}}
\put(4,-8){\line(0,1){6}}
\put(4,-8){\line(-1,2){8}}}
\begin{figure}[bottom]
	\begin{picture}(450,600)
\unitlength 10pt
\put(7.5,48){\lfig{\put(-7,0){\line(1,0){12}}
\put(-1,0){\line(0,1){6}}		\put(-1,0){\line(-2,-1){5.5}}
\put(5.5,-.25){$v_1$}	\put(-8,-.25){$v_3$}	\put(-1.5,6.5){$v_2$}
\put(-7,-3.25){$v_4$}	\put(3,4){$\sig_1$}	\put(-6,4){$\sig_2$}
\put(-8,-1.5){$\sig_3=C$}		\put(0,-4){$\sig_4$}}
{(a) The fan of the Hirzebruch surface with four top dimensional cones 
$\sig_i$ ($1\le i\le 4$). A chamber $C=\sig_3$ is chosen.}}

\put(22.5,48){\lfig{\hirzebruch
\put(-5.5,-3){$\xi_1(0,0)$}		\put(4.5,-2.25){$\xi_2(r,0)$}
\put(-5.25,8.5){$\xi_4(0,s)$}		\put(-1.5,-8.5){$\xi_3(r,s-ar)$}
\put(-4,-2){\circle*{.25}}		\put(-4,8){\circle*{.25}}
\put(4,-2){\circle*{.25}}		\put(4,-8){\circle*{.25}}
\thicklines
\put(-4,-2){\vector(1,0){2}}		\put(-4,-2){\vector(0,1){2}}
\put(4,-2){\vector(-1,0){2}}		\put(4,-2){\vector(0,1){2}}
\put(-4,8){\vector(0,-1){2}}		\put(-4,8){\vector(1,-2){2}}
\put(4,-8){\vector(0,-1){2}}		\put(4,-8){\vector(-1,2){2}}}
{(b) Position of weights $\xi_i$ in the fiber over $p_i\in F$.
The two arrows at $\xi_i$ are the negative of isotropy weights at $p_i$.}}

\put(37.5,48){\lfig{\hirzebruch
\put(-2,-7){$\bfgam{p_3,-C}$}		\put(1,5){$\bfgam{p_1,C}$}
\thicklines	\linethickness{1.5pt}
\put(-4.1,-2){\line(1,0){12}}		\put(-4,-2){\line(0,1){14}}
\put(4.05,-8){\line(-1,2){10}}		\put(3.95,-8){\line(-1,2){10}}
\put(4,-8){\line(0,-1){3}}
\multiput(-3.5,-1.5)(0,.5){17}{\circle*{.1}}
\multiput(-3,-1.5)(0,.5){15}{\circle*{.1}}
\multiput(-2.5,-1.5)(0,.5){13}{\circle*{.1}}
\multiput(-2,-1.5)(0,.5){11}{\circle*{.1}}
\multiput(-1.5,-1.5)(0,.5){9}{\circle*{.1}}
\multiput(-1,-1.5)(0,.5){7}{\circle*{.1}}
\multiput(-.5,-1.5)(0,.5){5}{\circle*{.1}}
\multiput(0,-1.5)(0,.5){3}{\circle*{.1}}
\put(.5,-1.5){\circle*{.1}}}
{(c) The two cones $\Gam^{p_1,C}$ and $\Gam^{p_3,-C}$, with 
$n^{p_1,C}=n^{p_3,-C}=0$, intersect at the dotted region.}}

\put(15,17){\lfig{\hirzebruch
\put(6,4){$\bfgam{p_2,C}$}		\put(-1,-6){$\bfgam{p_2,-C}$}
\put(-2,9){$\bfgam{p_4,C}$}		\put(-8,5){$\bfgam{p_4,-C}$}
\thicklines	\linethickness{1.5pt}
\put(4.06,-1.9){\line(1,0){5}}
\multiput(4.1,-1.9)(0,2){7}{\line(0,1){1}}
\put(3.9,-2.1){\line(0,-1){9}}
\multiput(3.97,-2.1)(-2.2,0){6}{\line(-1,0){1.1}}
\put(-3.9,8.02){\line(0,1){4}}
\multiput(-3.925,8)(1,-2){10}{\line(1,-2){.5}}
\multiput(-3.855,8.03)(1,-2){10}{\line(1,-2){.5}}
\put(-4.1,7.95){\line(-1,2){2}}		\put(-4.18,7.95){\line(-1,2){2}}
\multiput(-4.1,7.95)(0,-2.2){9}{\line(0,-1){1.1}}
\multiput(3.5,-2.5)(0,-.5){8}{\circle*{.1}}
\multiput(3,-2.5)(0,-.5){6}{\circle*{.1}}
\multiput(2.5,-2.5)(0,-.5){4}{\circle*{.1}}
\multiput(2,-2.5)(0,-.5){2}{\circle*{.1}}}
{(d) The four cones with polarizing indices $1$.
$\Gam^{p_4,C}$ and $\Gam^{p_2,-C}$ intersect at the dotted region.}}

\put(35,17){\lfig{\hirzebruch
\put(-8,-6){$\bfgam{p_1,-C}$}		\put(5,0){$\bfgam{p_3,C}$}
\thicklines	\linethickness{1.5pt}
\multiput(-4.05,-2.1)(-2,0){3}{\line(-1,0){1}}
\multiput(-4.1,-2.1)(0,-2){5}{\line(0,-1){1}}
\multiput(4.12,-8)(1,-2){2}{\line(1,-2){.5}}
\multiput(4.03,-8.04)(1,-2){2}{\line(1,-2){.5}}
\multiput(4.1,-8)(0,2){10}{\line(0,1){1}}}
{(e) The two cones $\Gam^{p_3,C}$ and $\Gam^{p_1,-C}$, with 
$n^{p_3,C}=n^{p_1,-C}=2$, have empty intersection.}}
	\end{picture}
\caption{\label{hirzebruch} Hirzebruch surface with an invariant line bundle,
case $s<ar$}
	\end{figure}

\setcounter{sect}{4}	\setcounter{figure}{0}

\begin{figure}[ht]
	\begin{picture}(425,360)
\unitlength 10pt
\put(10,24){\begin{picture}(20,20)
\put(0,10){\line(-1,-2){10}}		\put(0,10){\line(1,-2){10}}
\put(-10,-10){\line(1,0){20}}
\put(0,0){\line(-1,-2){2.5}}		\put(0,0){\line(1,-2){ 2.5}}
\put(-2.5,-5){\line(1,0){5}}		
\put(0,10){\line(-1,-6){2.5}}		\put(0,10){\line(0,-1){10}}	
\put(-10,-10){\line(5,2){12.5}}		\put(-10,-10){\line(3,2){7.5}}
\put(10,-10){\line(-3,2){7.5}}		\put(10,-10){\line(-1,1){10}}
\put(-.25,10.5){$a$}	\put(10.25,-10.5){$b$}	\put(-10.75,-10.5){$c$}
\put(.25,0){$d$}	\put(2.25,-5.75){$e$}	\put(-3.25,-5){$f$}
\put(-.25,-3.5){$\sig'$}
\put(-8,1){0}		\put(7,1){0}		\put(-.5,-13){0}
\put(-4.5,-3){1}	\put(1.75,1){2}		\put(2,-8){3}
\put(-1,2){4}		\put(3,-4.5){5}		\put(-2.5,-6.5){6}
\put(-9,-21){\parbox{180pt}{{\small 
(a) A fan that does not admit a strictly convex piecewise linear function.
The edges are spanned by $a(1,0,0)$, $b(0,1,0)$, $c(0,0,1)$,
$d(0,-2,-1)$, $e(-1,0,-2)$, $f(-2,-1,0)$.}}}
	\end{picture}}
\put(35,24){\begin{picture}(20,20)
\put(0,10){\line(-1,-2){10}}		\put(0,10){\line(1,-2){10}}
\put(-10,-10){\line(1,0){20}}
\put(-10,-10){\line(2,1){8}}		\put(-10,-10){\line(5,4){10}}
\put(-10,-10){\line(5,6){6.6667}}	\put(-10,-10){\line(2,3){8}}
\put(10,-10){\line(-3,1){12}}		\put(10,-10){\line(-2,1){8}}
\put(10,-10){\line(-5,4){10}}		\put(10,-10){\line(-3,4){6}}
\put(0,10){\line(-1,-4){2}}		\put(0,10){\line(0,-1){12}}
\put(0,10){\line(1,-6){1.3333}}		\put(0,10){\line(1,-3){4}}
\put(-2,2){\line(-1,-3){1.3333}}	\put(-2,2){\line(1,-2){4}}
\put(0,-2){\line(-1,-2){2}}		\put(0,-2){\line(1,3){1.3333}}
\put(4,-2){\line(-2,3){2.6667}}		\put(4,-2){\line(-1,0){7.3333}}
\put(-2,-6){\line(1,0){4}}
\put(-.25,10.5){$d$}	\put(10.25,-10.5){$e$}	\put(-10.75,-10.5){$f$}
\put(-1.75,2){$g$}	\put(3.65,-2.95){$h$}	\put(-2.1,-6.95){$i$}
\put(1.5,2){$j$}	\put(2.2,-5.8){$k$}	\put(-3.3,-2.95){$l$}
\put(.25,-1.75){$m$}		\put(-2,-.75){7}
\put(-9,-21){\parbox{185pt}{{\small (b) The triangulation of $\sig'$ 
by Jurkiewicz.
New edges are generated by $g(-1,-2,-1)$, $h(-1,-1,-2)$, $i(-2,-1,-1)$,
$j(-1,-2,-2)$, $k(-2,-1,-2)$, $l(-2,-2,-1)$, $m(-1,-1,-1)$.}}}
	\end{picture}}
	\end{picture}
\caption{\label{counter} Stereographic projection of the fan $\Sig$
corresponding to a non-\ka toric $3$-manifold.
There are 22 top-dimensional cones $\sig_i$ ($0\le i\le 21$) in $\Sig$.
The cones $\sig_i$ ($0\le i\le 6$) are shown in (a). 
The cones $\sig_i$ ($7\le i\le 21$) are contained in $\sig'$.
Only $\sig_7$ is labeled in (b).}
	\end{figure}
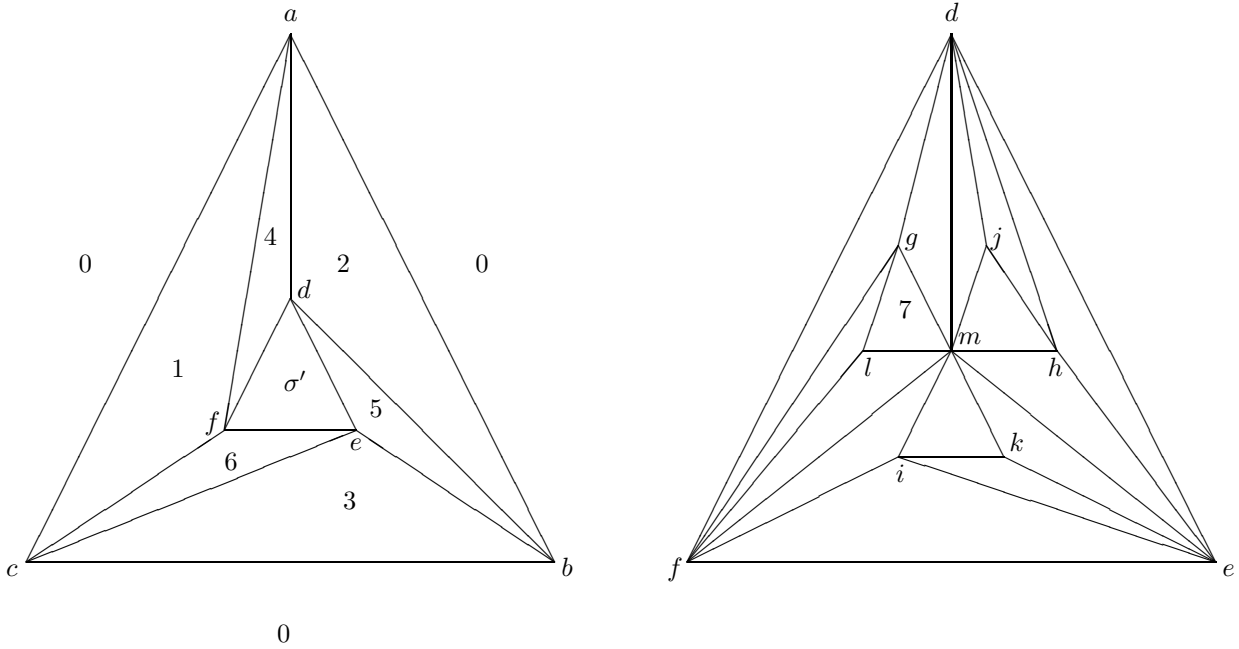

\newcommand{\tolman}{\put(-6,-3){\line(1,0){15}}
\put(-6,-3){\line(1,1){9}}		\put(-6,-3){\line(0,1){6}}
\put(3,3){\line(-1,0){9}}  \put(3,3){\line(0,1){3}}  \put(3,3){\line(1,-1){6}}
\put(-9,6){\line(1,0){12}}		\put(-9,6){\line(1,-1){3}}
\put(-9,6){\line(2,-1){18}}}
\begin{figure}[bottom]
	\begin{picture}(450,535)
\unitlength 10pt
\put(12.5,43.5){\fig{
\put(6,-3){\line(-1,0){12}}\put(6,-3){\line(-1,1){9}}\put(6,-3){\line(-2,1){6}}
\put(-3,0){\line(1,0){3}}\put(-3,0){\line(0,1){6}}\put(-3,0){\line(-1,-1){3}}
\put(-6,6){\line(1,0){3}}\put(-6,6){\line(0,-1){9}}\put(-6,6){\line(1,-1){6}}
\put(-6,-3){\circle*{.25}}\put(6,-3){\circle*{.25}}\put(-3,0){\circle*{.25}}
\put(0,0){\circle*{.25}}\put(-6,6){\circle*{.25}}\put(-3,6){\circle*{.25}}
\put(-6.5,-3.75){$\mu_1$}\put(-3,-.75){$\mu_2$}\put(-7,6.5){$\mu_3$}
\put(-3,6.5){$\mu_4$}\put(6,-3.75){$\mu_5$}\put(-.75,-.75){$\mu_6$}}
{(a) The moment polytope of Tolman's symplectic manifold.
The vertices $\mu_i$ are the values of the moment map $\mu$ at $p_i\in F$.}}

\put(32.5,43.5){\fig{\tolman
\put(-6,-3){\circle*{.25}}\put(3,6){\circle*{.25}}\put(-6,3){\circle*{.25}}
\put(3,3){\circle*{.25}}\put(9,-3){\circle*{.25}}\put(-9,6){\circle*{.25}}
\put(-7,-4){$\xi_1(0,0)$}		\put(-.25,6.5){$\xi_2(3,3)$}
\put(-9.25,2.25){$\xi_3(0,2)$}		\put(3.5,2.75){$\xi_4(3,3)$}
\put(8.3,-4){$\xi_5(5,0)$}		\put(-9,6.5){$\xi_6(-1,3)$}
\thicklines
\put(-6,-3){\vector(1,0){2}}		\put(-6,-3){\vector(1,1){2}}
\put(-6,-3){\vector(0,1){2}}
\put(-6,3){\vector(1,0){2}}		\put(-6,3){\vector(0,-1){2}}
\put(-6,3){\vector(1,-1){2}}
\put(-9,6){\vector(-1,0){2}}		\put(-9,6){\vector(-1,1){2}}
\put(-9,6){\vector(2,-1){4}}
\put(3,6){\vector(1,0){2}}		\put(3,6){\vector(0,1){2}}
\put(3,6){\vector(-1,-1){2}}
\put(3,3){\vector(-1,0){2}}		\put(3,3){\vector(0,-1){2}}
\put(3,3){\vector(1,-1){2}}
\put(9,-3){\vector(-1,0){2}}		\put(9,-3){\vector(-1,1){2}}
\put(9,-3){\vector(-2,1){4}}}
{(b) Position of the weights $\xi_i$ in the fiber of a line bundle 
over $p_i\in F$.
The negative of isotropy weights are also shown.}}

\put(12.5,15){\fig{\tolman
\put(-3,3){\circle*{.5}}		\put(-5.25,2){$(1,2)$}
\put(-8.5,0){$\bfgam{p_1,C}$}		\put(2,-2){$\bfgam{p_5,C}$}
\thicklines	\linethickness{1.5pt}
\put(9.05,-3){\line(-1,1){14}}		\put(8.96,-3.01){\line(-1,1){13.93}}
\put(9,-3){\line(-1,0){20}}		\put(-5.9,-2.77){\line(-1,0){5}}
\multiput(-5.85,-2.8)(2,2){6}{\line(1,1){1}}
\multiput(-5.9,-2.75)(2,2){6}{\line(1,1){1}}}
{(c) Two cones $\Gam^{p_1,C}$, $\Gam^{p_5,C}$ containing $e^*_1+2e^*_2$.
The polarizing indices $n^{p_1,C}=1$, $n^{p_5,C}=0$.}}

\put(32.5,15){\fig{\tolman
\put(-3,3){\circle*{.5}}		\put(-5.25,2){$(1,2)$}
\put(-4.5,4){$\bfgam{p_3,-C}$}		\put(-4,-2){$\bfgam{p_6,-C}$}
\thicklines	\linethickness{1.5pt}
\put(-6.1,3){\line(1,0){15}}		\put(-6,3){\line(0,-1){10.5}}
\multiput(-8.86,5.9)(2.4,0){8}{\line(1,0){1.2}}
\multiput(-8.75,5.91)(2,-2){7}{\line(1,-1){1}}
\multiput(-8.8,5.86)(2,-2){7}{\line(1,-1){1}}}
{(d) Two cones $\Gam^{p_3,-C}$, $\Gam^{p_6,-C}$ containing $e^*_1+2e^*_2$.
The polarizing indices $n^{p_3,-C}=0$, $n^{p_6,-C}=2$.}}
	\end{picture}
\caption{\label{tolman} Moment polytope and a $T^2$-invariant line
bundle over Tolman's manifold}
	\end{figure}

	\end{document}